\documentclass[aos,MSNbibl,seceqn,dvips]{arximspdf}

\doi{10.1214/13-AOS1154} 
\volume{41}
\issue{5}
\pubyear{2013}
\firstpage{2572}
\lastpage{2607}

\makeatletter
\newcommand{\rrVert}{\Vert}
\newcommand{\llVert}{\Vert}
\newtheorem{theorem}{Theorem}[section]
\newtheorem{lem}[theorem]{Lemma}
\newproclaim{rem}[theorem]{Remark}
\makeatother

\begin{document}
\begin{frontmatter}

\title{Convergence rates of eigenvector empirical spectral distribution
of large dimensional sample~covariance matrix}
\runtitle{Asymptotics of eigenvectors and eigenvalues}

\begin{aug}
\author[A]{\fnms{Ningning} \snm{Xia}\ead[label=e1]{xiann664@gmail.com}\thanksref{t1}},
\author[B]{\fnms{Yingli} \snm{Qin}\ead[label=e2]{yingli.qin@uwaterloo.ca}\thanksref{t2}}
\and
\author[A]{\fnms{Zhidong} \snm{Bai}\ead[label=e3]{baizd@nenu.edu.cn}\thanksref{t3}\corref{}}
\runauthor{N. Xia, Y. Qin and Z. D. Bai}
\affiliation{Northeast Normal University and National University of Singapore,\\
University of Waterloo, and
Northeast Normal University and\\ National University of Singapore}
\address[A]{N. Xia\\
Z. D. Bai\\
KLAS and School of Mathematics\\
\quad and Statistics\\
Northeast Normal University\\
Changchun 130024\\
China\\
and
Department of Statistics\\
\quad and Applied Probability\\
National University of Singapore\\
Singapore 117546\\
Singapore\\
\printead{e1}\\
\phantom{E-mail:\ }\printead*{e3}}
\address[B]{Y. Qin\\
Department of Statistics\\
\quad and Actuarial Science\\
University of Waterloo\\
Waterloo, ON N2L 3G1\\
Canada\\
\printead{e2}}
\end{aug}
\thankstext{t1}{Supported by the Fundamental Research Funds for the
Central Universities 11SSXT131.}
\thankstext{t2}{Supported by University of Waterloo Start-up Grant.}
\thankstext{t3}{Supported by NSF China Grant 11171057, PCSIRT Fundamental Research Funds for the Central Universities,
and NUS Grant R-155-000-141-112.}

\received{\smonth{2} \syear{2013}}
\revised{\smonth{6} \syear{2013}}

%
\begin{abstract}
The eigenvector Empirical Spectral Distribution (VESD) is adopted to
investigate the limiting behavior of eigenvectors and eigenvalues of
covariance matrices. In this paper, we shall show that the Kolmogorov
distance between the expected VESD of sample covariance matrix and the
Mar\v{c}enko--Pastur distribution function is of order $O(N^{-1/2})$.
Given that data dimension $n$ to sample size $N$ ratio is bounded
between 0 and 1, this convergence rate is established under finite 10th
moment condition of the underlying distribution. It is also shown that,
for any fixed $\eta>0$, the convergence rates of VESD are $O(N^{-1/4})$
in probability and $O(N^{-1/4+\eta})$ almost surely, requiring finite
8th moment of the underlying distribution.
\end{abstract}

%

\end{frontmatter}

\section{Introduction and main results}\label{sec1}

Let $\mathbf{X}_i=(X_{1i},X_{2i},\ldots, X_{ni})^T$ and $\mathbf
{X}=(\mathbf{X}_1,\ldots,\mathbf{X}_N)$ be an $n\times N$ matrix of
i.i.d. (independent and identically distributed) complex random
variables with mean 0 and variance 1. We consider, a class of sample
covariance matrices 
\[
\mathbf{S}_n=\frac{1}{N}\sum_{k=1}^N
\mathbf{X}_k\mathbf{X}_k^*=\frac{1}{N} \mathbf{X}
\mathbf{X}^*,
\]
where $\mathbf{X}^*$ denotes the conjugate transpose
of the data matrix $\mathbf{X}$. The Empirical
Spectral Distribution (ESD) $F^{\mathbf{S}_n}(x)$ of $\mathbf{S}_n$ is
then defined as
%
%
\begin{equation}\label{esd}
F^{\mathbf{S}_n}(x)=\frac{1}{n}\sum_{i=1}^n \mathrm{I}(\lambda_i \leq
x),
\end{equation}
where $\lambda_1\leq\cdots\leq\lambda_n$ are the eigenvalues of
$\mathbf{S}_n$ in ascending order and $\mathrm{I}(\cdot)$ is the
conventional indicator function.

Mar\v{c}enko and Pastur in \cite{mplaw1967} proved that with
probability 1, $F^{\mathbf{S}_n}(x)$ converges weakly to the standard
Mar\v{c}enko--Pastur distribution $F_y(x)$ with density function
%
%
\begin{equation}\label{mplaw}
p_y(x)=\frac{dF_y(x)}{dx}=\frac{1}{2\pi xy}\sqrt{(x-a) (b-x)}
\mathrm{I}(a\leq x\leq b),
\end{equation}
where $a=(1-\sqrt{y})^2$ and $b=(1+\sqrt{y})^2$. Here the positive
constant $y$ is the limit of dimension to sample size ratio when both
$n$ and $N$ tend to infinity.

In applications of asymptotic theorems of spectral analysis of large
dimensional random matrices, one of the important problems is the
convergence rate of the ESD. The Kolmogorov distance between the
expected ESD of $\mathbf{S}_n$ and the Mar\v{c}enko--Pastur
distribution $F_y(x)$ is defined as
\[
\Delta=\bigl\|\mathrm{E}F^{\mathbf{S}_n}-F_y\bigr\|=\sup_x\bigl|
\mathrm{E}F^{\mathbf{S}_n}(x)-F_y(x)\bigr|
\]
as well as the distance between two distributions $F^{\mathbf{S}_n}(x)$ and $F_y(x)$,
\[
\Delta_p=\bigl\|F^{\mathbf{S}_n}-F_y\bigr\|=\sup
_x\bigl|F^{\mathbf{S}_n}(x)-F_y(x)\bigr|.
\]
Notice that, for any constant $C>0$,
\[
P(\Delta_p\geq C)=P \Bigl\{\sup_x\bigl|F^{\mathbf{S}_n}(x)-F_y(x)\bigr|
\geq C \Bigr\}\leq C^{-1}\mathrm{E}\Delta_p.
\]
Thus, $\Delta_p$ measures the rate of convergence in probability.

Bai in \cite{bai1993a,bai1993b} firstly tackled the problem of
convergence rate and established three Berry--Esseen type inequalities
for the difference of two distributions in terms of their Stieltjes
transforms. G\"{o}tze and Tikhomirov in \cite{gt2004} further improved
the Berry--Esseen type inequlatiy and showed the convergence rate of
$F^{\mathbf{S}_n}(x)$ is $O(N^{-1/2})$ in probability under finite 8th moment
condition. More recently, a sharper bound is obtained by Pillai and Yin
in \cite{yinev}, under a stronger condition, that is, the
sub-exponential decay assumption. It is shown that the difference
between eigenvalues of $\mathbf{S}_n$ and the Mar\v{c}enko--Pastur
distribution is of order $O (N^{-1}(\log N)^{O(\log\log N)} )$ in
probability.

In the literature, research on limiting properties of eigenvectors of
large dimensional sample covariance matrices is much less developed
than that of eigenvalues, due to the cumbersome formulation of the
eigenvectors. Some great achievements have been\vadjust{\goodbreak} made in proving the
properties of eigenvectors for large dimensional sample covariance
matrices, such as \cite{bmp2007,jack1981,jack1984,jack1989,jack1990},
and that for Wigner matrices, such as \cite{je,te,erdos}.

However, the eigenvectors of large sample covariance matrices play an
important role in high-dimensional statistical analysis. In particular,
due to the increasing availability of high-dimensional data, principal
component analysis (PCA) has been favorably recognized as a powerful
technique to reduce dimensionality. The eigenvectors corresponding to
the leading eigenvalues are the directions of the principal components.
Johnstone \cite{ij2001} proposed the spiked eigenvalue model to test
the existence of principal component. Paul \cite{paul} discussed the
length of the eigenvector corresponding to the spiked eigenvalue.

In PCA, the eigenvectors ($\bolds{\nu}_1^0,\ldots,\bolds{\nu}_n^0$) of
population covariance matrix $\bolds{\Sigma}$ determine the directions
in which we project the observed data and the corresponding eigenvalues
($\lambda_1^0,\ldots,\lambda_n^0$) determine the proportion of total
variability loaded on each direction of projections. In practice, the
(sample) eigenvalues ($\lambda_1,\ldots,\lambda_n$) and eigenvectors
($\bolds{\nu} _1,\ldots,\bolds{\nu}_n$) of the sample covariance
matrix $\mathbf{S}_n$ are used in PCA. In \cite{anderson1963}, Anderson
has shown the following asymptotic distribution for the sample
eigenvectors $\bolds{\nu}_1,\ldots,\bolds{\nu}_n$ when the observations
are from a multivariate normal distribution of covariance matrix
$\bolds{\Sigma}$ with distinct eigenvalues:
\[
\sqrt{N}\bigl(\bolds{\nu}_i-\bolds{\nu}_i^0\bigr)
\stackrel{d} {\to}N_n(0,\mathbf{D}_i),
\]
where
\[
\mathbf{D}_i=\lambda_i^0\sum
_{k=1,k\neq i}^n\frac{\lambda_k^0}{(\lambda
_k^0-\lambda_i^0)^2}\bolds{\nu}_k^0{
\bolds{\nu}_k^0}^T.
\]

However, this is a large sample result when the dimension $n$ is fixed
and low. In particular, if $\bolds{\Sigma}=\sigma^2\mathbf{I}_n$, then
the eigenmatrix (matrix of eigenvectors) should be asymptotically
isotropic when the sample size is large. That is, the eigenmatrix
should be asymptotically Haar, under some minor moment conditions.
However, when the dimension is large (increasing), the Haar property is
not easy to formulate.

Motivated by the orthogonal iteration method, \cite{ma2013} proposed an
iterative thresholding method to estimate sparse principal subspaces
(spanned by the leading eigenvectors of $\bolds{\Sigma}$) in high
dimensional and spiked covariance matrix setting. The convergence rates
of the proposed estimators are provided. By reducing the sparse PCA
problem to a high-dimensional regression problem, \cite{cai2013}
established the optimal rates of convergence for estimating the
principal subspace with respect to a large collection of spiked
covariance matrices. See the reference therein for more literature on
sparse PCA and spiked covariance matrices.

To perform the test of existence of spiked eigenvalues, one has to
investigate the null properties\vadjust{\goodbreak} of the eigenmatrices, that is, when
$\bolds{\Sigma}=\sigma^2\mathbf{I}_n$ (i.e., nonspiked). Then the
eigenmatrix should be asymptotically isotropic, when the sample size is
large. That is, the eigenmatrix should be asymptotically Haar. However,
when the dimension is large, the Haar property is not easy to
formulate. The recent development in random matrix theory can help us
investigate the large dimension and large sample properties of
eigenvectors. We will adopt the VESD, defined later in the paper, to
characterize the asymptotical Haar property so that if the eigenmatrix
is Haar, then the process defined the VESD tends to a Brownian bridge.
Conversely, if the process defined by the VESD tends to a Brownian
bridge, then it indicates a similarity between the Haar distribution
and that of the eigenmatrix. Therefore, studying the large sample and
large dimensional results of the VESD can assist us in better examining
spiked covariance matrix as assumed by \cite{ma2013} and \cite{cai2013}
among many others.


Let $\mathbf{U}_n\bolds{\Lambda}_n\mathbf{U}_n^*$ denote the spectral
decomposition of $\mathbf{S}_n$, where
$\bolds{\Lambda}_n=\break\operatorname{diag}(\lambda_1,\lambda_2,
\ldots,\lambda_n)$ and $\mathbf{U}_n=(u_{ij})_{n\times n}$ is a unitary
matrix consisting of the corresponding orthonormal eigenvectors of
$\mathbf{S}_n$. For each $n$, let $\mathbf{x}_n\in\mathbb{C}^n$,
$\|\mathbf{x}_n\|=1$ be nonrandom and let
$\mathbf{d}_n=\mathbf{U}_n^*\mathbf{x}_n=(d_1,\ldots,d_n)^*$, where
$\|\mathbf{x}_n\|$ denotes Euclidean norm of $\mathbf{x}_n$.

Define a stochastic process $X_n(t)$ by
\[
X_n(t)=\sqrt{n/2}\sum_{j=1}^{[nt]}
\biggl(|d_j|^2-\frac1n \biggr),\qquad[a]\mbox{ denotes the greatest integer}\leq a.
\]
If $\mathbf{U}_n$ is Haar distributed over the orthogonal matrices,
then $\mathbf{d}_n$ would be uniformly distributed over the unit sphere
in $\mathbb{R}^n$, and the limiting distribution of $X_n(t)$ is a
unique Brownian bridge $B(t)$ when $n$ tends to infinity.
In this paper, we use the behavior of $X_n(t)$ for all $\mathbf{x}_n$
to reflect the uniformity of $\mathbf{U}_n$. The process $X_n(t)$ is
considerably important for us to understand the behavior of the
eigenvectors of $\mathbf{S}_n$.

Motivated by Silverstein's ideas in
\cite{jack1981,jack1984,jack1989,jack1990}, we want to examine the
limiting properties of $\mathbf{U}_n$ through stochastic process
$X_n(t)$. We claim that $\mathbf{U}_n$ is ``asymptotically Haar
distributed,'' which means $X_n(t)$ converges to a Brownian bridge
$B(t)$. In \cite{jack1989}, it showed that the weak convergence of
$X_n(t)$ converging to a Brownian bridge $B(t)$ is equivalent to $X_n
(F^{\mathbf{S}_n}(x) )$ converging to $B (F_{y}(x) )$. We therefore consider
transforming $X_n(t)$ to $ X_n (F^{\mathbf{S}_n}(x) ) $ where $F^{\mathbf{S}_n}(x)$ is
the ESD of $\mathbf{S}_n$.

We define the eigenvector Empirical Spectral Distribution (VESD)
$H^{\mathbf{S}_n}(x)$ of $\mathbf{S}_n$ as follows:
%
%
\begin{equation}
\label{esed} H^{\mathbf{S}_n}(x)=\sum_{i=1}^n
|d_i|^2\mathrm{I}(\lambda_i\leq x).
\end{equation}
Between $H^{\mathbf{S}_n}(x)$ in (\ref{esed}) and $F^{\mathbf{S}_n}(x)$
in (\ref{esd}), we notice that there is no difference except the
coefficient associated with each indicator function such that
%
%
\begin{equation}\label{haarconj2}
X_n\bigl(F^{\mathbf{S}_n}(x)\bigr)=\sqrt{n/2}\bigl(H^{\mathbf
{S}_n}(x)-F^{\mathbf{S}_n}(x)
\bigr).
\end{equation}
Henceforth, the investigation of $X_n(t)$ is converted to that of the
difference between two empirical distributions $H^{\mathbf{S}_n}(x)$
and $F^{\mathbf{S}_n}(x)$. The authors in \cite{bmp2007} proved that
$H^{\mathbf{S}_n}(x)$ and $F^{\mathbf{S}_n}(x)$ have the same limiting
distribution, the Mar\v{c}enko--Pastur distribution $F_y(x)$, where
$y_n=n/N$ and $y=\lim_{n, N\to\infty}y_n\in(0,1)$.

Before we present the main theorems, let us introduce the following
notation:
\[
\Delta^H=\bigl\|\mathrm{E}H^{\mathbf{S}_n}-F_{y_n}\bigr\|=\sup
_x\bigl|\mathrm{E}H^{\mathbf{S}_n}(x)-F_{y_n}(x)\bigr|
\]
and
\[
\Delta^H_p=\bigl\|H^{\mathbf{S}_n}-F_{y_n}\bigr\|=\sup
_x\bigl|H^{\mathbf{S}_n}(x)-F_{y_n}(x)\bigr|.
\]
We denote $\xi_n=O_p(a_n)$ and $\eta_n=O_{\mathrm{a.s.}}(b_n)$ if, for
any $\epsilon>0$, there exist a large positive constant $c_1$ and a
positive random variable $c_2$, such that
\[
P (\xi_n/a_n\geq c_1 )\leq\epsilon\quad
\mbox{and}\quad P(\eta_n/b_n\leq c_2)=1,
\]
respectively.

In this paper, we follow the work in \cite{bmp2007} and establish three
types of convergence rates of $H^{\mathbf{S}_n}(x)$ to $F_{y_n}(x)$ in the
following theorems.

%
%
\begin{theorem}\label{theorem1}
Suppose that $X_{ij}$, $i=1,\ldots, n$, $j=1,\ldots,N$ are i.i.d.
complex random variables with $\mathrm{E}X_{11}=0$,
$\mathrm{E}|X_{11}|^2=1$ and $\mathrm{E}|X_{11}|^{10}< \infty$. For any
fixed unit vector $\mathbf{x}_n\in
\mathbb{C}_1^n=\{\mathbf{x}\in\mathbb{C}^n\dvtx \|\mathbf{x}\|=1\}$,
and $y_n=n/N\leq1$, it then follows that
\[
\Delta^H=\bigl\|\mathrm{E}H^{\mathbf{S}_n}-F_{y_n} \bigr\|=
\cases{ O\bigl(N^{-1/2}a^{-3/4}\bigr), &\quad if $N^{-1/2}
\leq a<1$, \vspace*{2pt}
\cr
O\bigl(N^{-1/8}\bigr),&\quad if
$a<N^{-1/2}$,}
\]
where $a=(1-\sqrt{y_n})^2$ as it is defined in (\ref{mplaw}) and
$F_{y_n}$ denotes the Mar\v{c}enko--Pastur distribution function with
an index $y_n$.
\end{theorem}

%
%
\begin{rem}
From the proof of Theorem \ref{theorem1}, it is clear that the
condition $\mathrm{E}|X_{11}|^{10}< \infty$ is required only in the
truncation step in the next section. We therefore believe that the
condition $\mathrm{E}|X_{11}|^{10}< \infty$ can be replaced by
$\mathrm{E}|X_{11}|^{8}< \infty$ in Theorems
\ref{theorem2}~and~\ref{theorem3}.
\end{rem}

%
\begin{rem}
Because the convergence rate of $\|\mathrm{E}H^{\mathbf{S}_n}-F_{y} \|
$ depends on the convergence rate of $|y_n-y|$, we only consider the
convergence rate of \mbox{$\|\mathrm{E}H^{\mathbf{S}_n}-F_{y_n} \|$}.
\end{rem}

%
%
\begin{rem}
As $a=(1-\sqrt{y_n})^2$, we can characterize the closeness between
$y_n$ and 1 through $a$.\vadjust{\goodbreak} In particular, when $y_n$ is away from 1 (or
$a\geq N^{-1/2}$), the convergence rate of
$\|\mathrm{E}H^{\mathbf{S}_n}-F_{y_n}\|$ is $O(N^{-1/2})$, which we believe is
the optimal convergence rate. This is because we observe in
\cite{bmp2007} that for an analytic function~$f$,
%
\begin{equation}\label{eq001}
Y_n(f)=\sqrt{n}\int f(x)\,d \bigl(H^{\mathbf{S}_n}(x)-F_{y_n}(x)
\bigr)
\end{equation}
converges to a Gaussian distribution. While in \cite{baijack2004}, Bai
and Silverstein proved that the limiting distribution of
\[
n\int f(x)\,d \bigl(F^{\mathbf{S}_n}(x)-F_{y_n}(x) \bigr)
\]
is also a Gaussian distribution.
We therefore conjecture that the optimal rate of $H^{\mathbf{S}_n}(x)$ should be
$O(N^{-1/2})$ and $O(N^{-1})$ for $F^{\mathbf{S}_n}(x)$. Although $F^{\mathbf{S}_n}(x)$
and $H^{\mathbf{S}_n}(x)$ converge to the same limiting distribution, there
exists a substantial difference between $F^{\mathbf{S}_n}(x)$ and $H^{\mathbf{S}_n}(x)$.
\end{rem}

%
%
\begin{rem}
Notice that two matrices $ \mathbf{X}\mathbf{X}^*$ and
$\mathbf{X}^*\mathbf{X}$ share the same set of nonzero eigenvalues.
However, these two matrices do not always share the same set of
eigenvectors.
Especially when $y_n\gg1$, the eigenvectors of $\mathbf{S}_n$
corresponding to 0 eigenvalues can be arbitrary. As a result, the limit
of $H^{\mathbf{S}_n}$ may not exist or heavily depends on the choice of unit
vector $\mathbf{x}_n$. Therefore, we only consider the case of $y_n\leq
1$ in this paper and leave the case of $y_n\geq1$ as a future research
problem.
\end{rem}

The rates of convergence in probability and almost sure convergence of
the VESD are provided in the next two theorems.

%
%
\begin{theorem}\label{theorem2}
Under the assumptions in Theorem \ref{theorem1} except that we now only
require $\mathrm{E}|X_{11}|^8< \infty$, we have
\[
\Delta_p^H=\bigl\|H^{\mathbf{S}_n}-F_{y_n} \bigr\|=
\cases{ O_{p}\bigl(N^{-1/4}a^{-1/2}\bigr), &\quad if
$N^{-1/4}\leq a<1$,
\vspace*{2pt}\cr
O_{p}
\bigl(N^{-1/8}\bigr), &\quad if $a<N^{-1/4}$.}
\]
\end{theorem}

%
%
\begin{rem}As an application of Theorem \ref{theorem2}, in \cite
{baixia2013} we extended the CLT of the linear spectral statistics
$ Y_n(f)$ established in \cite{bmp2007} to the case where the kernel
function $f$ is continuously twice differentiable provided that the
sample covariance matrix $\mathbf{S}_n$ satisfies the assumptions of
Theorem \ref{theorem2}. This result is useful in testing Johnstone's
hypothesis when
normality is not assumed.
\end{rem}

%
%
\begin{theorem}\label{theorem3}
Under the assumptions in Theorem \ref{theorem2}, for any
$\eta>0$, we have
\[
\Delta_p^H=\bigl\|H^{\mathbf{S}_n}-F_{y_n} \bigr\|=
\cases{O_{\mathit{a.s.}}\bigl(N^{-1/4+\eta}a^{-1/2}\bigr), &\quad if
$N^{-1/4}\leq a<1$, \vspace*{2pt}
\cr
O_{\mathit{a.s.}} \bigl(N^{-1/8+\eta}\bigr), &\quad if $a<N^{-1/4}$.}
\]
\end{theorem}

%
%
\begin{rem}
In this paper, we will use the following notation:
\begin{itemize}
\item $\mathbf{X}^*$ denote the conjugate transpose of a matrix (or
    vector) $\mathbf{X}$;

\item $\mathbf{X}^T$ denote the (ordinary) transpose of a matrix (or
    vector) $\mathbf{X}$;

\item $\|\mathbf{x}\|$ denote the Euclidean norm for any vector
    $\mathbf{x}$;

\item $\|\mathbf{A}\|=\sqrt{\lambda_{\max}(\mathbf{AA}^*)}$, the
    spectral norm;

\item $\|F\|=\sup_x |F(x)|$ for any function $F$;

\item $\bar{z}$ denote the conjugate of a complex number $z$.
\end{itemize}
\end{rem}

The rest of the paper is organized as follows. In Section~\ref{sec2},
we introduce the main tools used to prove Theorems \ref{theorem1},
\ref{theorem2} and \ref{theorem3}, including Stieltjes transform and a
Berry--Esseen type inequality. The proofs of these three theorems are
presented in Sections~\ref{sec3}--\ref{sec6}. Several important results
which are repeatedly employed throughout
Sections~\ref{sec3}--\ref{sec6} are proved in Appendix~\ref{sec7}.
Appendix~\ref{sec8} contains some existing results in the literature.
Finally, preliminaries on truncation, centralization and rescaling are
postponed to the last section.

\section{Main tools}\label{sec2}
\subsection{Stieltjes transform}\label{sec2.1}

The Stieltjes transform is an essential tool in random matrix theory
and our paper. Let us now briefly review the Stieltjes transform and
some important and relevant results. For a cumulative distribution
function $G(x)$, its Stieltjes transform $m_G(z)$ is defined as
\[
m_G(z)=\int\frac{1}{\lambda-z}\,dG(\lambda),\qquad z\in\mathbb{C}^+=
\bigl\{z\in\mathbb{C}, \Im(z)>0\bigr\},
\]
where $\Im(\cdot)$ denotes the imaginary part of a complex number. The
Stieltjes transforms of the ESD $F^{\mathbf{S}_n}(x)$ and the VESD
$H^{\mathbf{S}_n}(x)$ are
\[
m_{F^{\mathbf{S}_n}}(z)=\frac{1}{n}\operatorname{tr} (\mathbf{S}_n-z
\mathbf{I}_n)^{-1}
\]
and
\[
m_{H^{\mathbf{S}_n}}(z)=\mathbf{x}_n^* (\mathbf{S}_n-z
\mathbf{I}_n)^{-1}\mathbf{x}_n,
\]
respectively. Here $\mathbf{I}_n$ denotes the $n\times n$ identity
matrix. For simplicity of notation, we use $m_n(z)$ and $m_n^H(z)$ to
denote $ m_{F^{\mathbf{S}_n}}(z)$ and $m_{H^{\mathbf{S}_n}}(z)$, respectively.

%
%
\begin{rem}
Notice that although the eigenmatrix $\mathbf{U}_n$ may not be unique,
the Stieltjes transform $m_n^H(z)$ of $H^{\mathbf{S}_n}$ depends on
$\mathbf{S}_n$ for any $\mathbf{x}_n$ rather than $\mathbf{U}_n$.
\end{rem}

Let $\underline{\mathbf{S}}_n=\mathbf{X}^*\mathbf{X}/N$ denote the
companion matrix of $\mathbf{S}_n$. As $\mathbf{S}_n$ and
$\underline{\mathbf{S}}_n$ share the same set of nonzero eigenvalues,
it can be shown that Stieltjes transforms of $F^{\mathbf{S}_n}(x)$ and
$F^{\underline{\mathbf{S}}_n}(x)$ satisfy the following equality:
%
%
\begin{equation}\label{mnmn}
\underline{m}_n(z)=-\frac{1-y_n}{z}+y_n m_n(z),
\end{equation}
where $\underline{m}_n(z)$ denotes the Stieltjes transform of
$F^{\underline{\mathbf{S}}_n}(x)$. Moreover, \cite{bai1998} and
\cite{jackbai1995} claimed that $F^{\underline{\mathbf{S}}_n}$
converges, almost surely, to a nonrandom distribution function
$\underline{F}_y(x)$ with Stieltjes transform $\underline{m}(z)$ such
that
%
%
\begin{equation}\label{mm}
\underline{m}(z)=-\frac{1-y}{z}+y m_y(z),
\end{equation}
where $m_y(z)$ denotes the Stieltjes transform of the
Mar\v{c}enko--Pastur distribution with index $y$.
Using (6.1.4) in \cite{baijack2010}, we also obtain the relationship
between two limits $m_y(z)$ and $\underline{m}(z)$ as follows:
%
%
\begin{equation}\label{underlinemm}
m_y(z)=\frac{1}{-z
(1+\underline{m}(z) )}.
\end{equation}

\subsection{A Berry--Esseen type inequality}\label{sec2.2}
%
%
\begin{lem}\label{main}
Let $H^{\mathbf{S}_n}(x)$ and $F_{y_n}(x)$ be the VESD of $\mathbf{S}_n$ and the
Mar\v{c}enko--Pastur distribution with index $y_n$, respectively.
Denote their corresponding Stieltjes transforms by $m_n^H(z)$ and
$m_{y_n}(z)$, respectively. Then there exist large positive constants
$A,B,K_1,K_2$ and $K_3$, such that for $A>B>5$,
\begin{eqnarray*}
\Delta^H&=&\bigl\|\mathrm{E}H^{\mathbf{S}_n}(x)-F_{y_n}(x)\bigr\|
\\
&\leq& K_1\int_{-A}^A\bigl|
\mathrm{E}m_n^H(z)-m_{y_n}(z)\bigr|\,du
+K_2v^{-1}\int_{|x|>B}\bigl|
\mathrm{E}H^{\mathbf{S}_n}(x)-F_{y_n}(x)\bigr|\,dx
\\
&&{}+K_3v^{-1}\sup_x\int
_{|t|<v}\bigl|F_{y_n}(x+t)-F_{y_n}(x)\bigr|\,dt,
\end{eqnarray*}
where $z=u+iv$ is a complex number with positive imaginary part (i.e., $v>0$).
\end{lem}

%
%
\begin{rem}
Lemma \ref{main} can be proved using Lemma \ref{LL}. To prove Theorem
\ref{theorem1}, we apply Lemma \ref{main}. In addition, we prove
Theorems \ref{theorem2} and \ref{theorem3} by replacing
$\mathrm{E}H^{\mathbf{S}_n}(x)$, $\mathrm{E}m_n^H(z)$ with $H^{\mathbf{S}_n}(x)$ and
$m_n^H(z)$, respectively.
\end{rem}

\section{\texorpdfstring{Proof of Theorem \protect\ref{theorem1}}{Proof of Theorem 1.1}}\label{sec3}

Under the condition of $\mathrm{E}|X_{11}|^{10} < \infty$, we can
choose a sequence of $\eta_N$ with $\eta_N \downarrow0$ and
$\eta_NN^{1/4}\uparrow\infty$ as $N\to\infty$, such that
%
%
\begin{equation}\label{eta}
\lim_{N \rightarrow\infty} \frac{1}{\eta_N^{10}} \mathrm{E} \bigl
(|X_{11}|^{10}
\mathrm{I}\bigl(|X_{11}| > \eta_NN^{1/4}\bigr)
\bigr)=0.
\end{equation}
Furthermore, without loss of generality, we can assume that every
$|X_{ij}|$ is bounded by $\eta_NN^{1/4}$ and has mean 0 and variance 1.
See Appendix~\ref{sec9} for details on truncation, centralization and
rescaling.

We introduce some notation before start proving Theorem
\ref{theorem1}. Throughout the paper, we use $C$ and $C_i$ for
$i=0,1,2,\ldots$ to denote positive constant numbers which are
independent of $N$ and may take different values at different
appearances. Let $\mathbf{X}_j$ denote the $j$th column of the data
matrix $\mathbf{X}$. Let $\mathbf{r}_j= \mathbf{X}_j/\sqrt{N}$ so that
$ \mathbf{S}_n=\sum_{j=1}^N \mathbf{r}_j\mathbf{r}_j^*$ and let
\begin{eqnarray*}
v_y&=&\sqrt{a}+\sqrt{v}=1-\sqrt{y_n}+\sqrt{v},
\\
\mathbf{B}_j&=&\mathbf{S}_n-\mathbf{r}_j\mathbf{r}_j^*,
\\
\mathbf{A}(z)&=&\mathbf{S}_n-z\mathbf{I}_n,
\\
\mathbf{A}_j(z)&=&\mathbf{B}_j-z\mathbf{I}_n,
\\
\alpha_j(z)&=&\mathbf{r}_j^*\mathbf{A}_j^{-1}(z)\mathbf{x}_n\mathbf{x}_n^*\mathbf{r}_j -\frac{1}{N}\mathbf{x}_n^*\mathbf{A}_j^{-1}(z)\mathbf{x}_n,
\\
\xi_j(z)&=&\mathbf{r}_j^*\mathbf{A}_j^{-1}(z)\mathbf{r}_j-\frac{1}{N}\mathrm{E}\operatorname{tr}\mathbf{A}_j^{-1}(z),
\\
\hat{\xi}_j(z)&=&\mathbf{r}_j^*\mathbf{A}_j^{-1}(z)\mathbf{r}_j-\frac{1}{N}\operatorname{tr}\mathbf{A}_j^{-1}(z),
\\
b(z)&=&\frac{1}{1+(1/N)\mathrm{E}\operatorname{tr}\mathbf{A}^{-1}(z)},
\\
b_1(z)&=&\frac{1}{1+(1/N)\mathrm{E}\operatorname{tr}\mathbf{A}_1^{-1}(z)},
\\
\beta_j(z)&=&\frac{1}{1+\mathbf{r}_j^*\mathbf{A}_j^{-1}(z)\mathbf{r}_j}.
\end{eqnarray*}

It is easy to show that
%
%
\begin{equation}\label{betabxi}
\beta_j(z)-b_1(z)=-b_1(z)\beta_j(z)\xi_j(z).
\end{equation}
For any $j=1,2,\ldots,N$, we can also show that
%
%
\begin{equation}\label{raaj}
\mathbf{r}_j^*\mathbf{A}^{-1}(z)=\beta_j(z)\mathbf{r}_j^*\mathbf{A}_j^{-1}(z)
\end{equation}
due to the fact that
\[
(\mathbf{B}_j-z\mathbf{I}_n)^{-1}-\bigl(
\mathbf{B}_j+\mathbf{r}_j\mathbf{r}_j^*-z
\mathbf{I}_n\bigr)^{-1} =(\mathbf{B}_j-z
\mathbf{I}_n)^{-1}\mathbf{r}_j
\mathbf{r}_j^*\bigl(\mathbf{B}_j+
\mathbf{r}_j\mathbf{r}_j^*-z\mathbf{I}_n
\bigr)^{-1}.
\]
From (2.2) in \cite{jack1995}, we can write $\underline{m}_n(z)$ in
terms of $\beta_j(z)$ as follows:
%
%
\begin{eqnarray}\label{mbeta}
\underline{m}_n(z)&=&-\frac{1}{zN}\sum
_{j=1}^N\beta_j(z).
\end{eqnarray}
We proceed with the proof of Theorem \ref{theorem1}:
\begin{eqnarray*}
\delta&=:&\mathrm{E}m_n^H(z)-m_y(z)
\\
&=&\mathbf{x}_n^* \bigl[\mathrm{E}\mathbf{A}^{-1}(z)-
\bigl(-z\underline{m}(z)\mathbf{I}_n-z\mathbf{I}_n
\bigr)^{-1} \bigr]\mathbf{x}_n
\\
&=& \bigl(z\underline{m}(z)+z \bigr)^{-1}\mathbf{x}_n^*
\mathrm{E} \bigl[ \bigl(z\underline{m}(z)+z \bigr)\mathbf{A}^{-1}(z)+
\mathbf{I}_n \bigr]\mathbf{x}_n
\\
&=& \bigl(z\underline{m}(z)+z \bigr)^{-1}\mathbf{x}_n^*
\mathrm{E} \bigl[ \bigl(z\mathbf{I}_n+\mathbf{A}(z) \bigr)
\mathbf{A}^{-1}(z)+z\underline{m}(z)\mathbf{A}^{-1}(z) \bigr]
\mathbf{x}_n
\\
&=& \bigl(z\underline{m}(z)+z \bigr)^{-1}\mathbf{x}_n^*
\mathrm{E} \Biggl[\sum_{j=1}^N
\mathbf{r}_j\mathbf{r}_j^*\mathbf{A}^{-1}(z)+z
\bigl(\mathrm{E}\underline{m}_n(z) \bigr)\mathbf{A}^{-1}(z)
\\
&&\hspace*{85pt}{} -z \bigl(\mathrm{E}\underline{m}_n(z) \bigr)
\mathbf{A}^{-1}(z)+z\underline{m}(z)\mathbf{A}^{-1}(z)
\Biggr]\mathbf{x}_n
\\
&=& \bigl(z\underline{m}(z)+z \bigr)^{-1}\mathbf{x}_n^*
\mathrm{E} \Biggl[\sum_{j=1}^N
\beta_j(z)\mathbf{r}_j\mathbf{r}_j^*
\mathbf{A}_j^{-1}(z)- \bigl(-z\mathrm{E}
\underline{m}_n(z) \bigr)\mathbf{A}^{-1}(z)
\\[-2pt]
&&\hspace*{150pt}{} - \bigl(z\mathrm{E}\underline{m}_n(z)-z\underline{m}(z)\bigr)\mathbf{A}^{-1}(z) \Biggr]\mathbf{x}_n
\\
&=& \bigl(z\underline{m}(z)+z \bigr)^{-1}\mathbf{x}_n^*
\mathrm{E} \Biggl[\sum_{j=1}^N
\beta_j(z)\mathbf{r}_j\mathbf{r}_j^*
\mathbf{A}_j^{-1}(z)- \Biggl(\frac{1}{N}\sum
_{j=1}^N\mathrm{E}\beta_j(z)
\Biggr)\mathbf{A}^{-1}(z) \Biggr]\mathbf{x}_n
\\
&&{}- \bigl(z\underline{m}(z)+z \bigr)^{-1}\mathbf{x}_n^*\bigl(z\mathrm{E}\underline{m}_n(z)-z\underline{m}(z) \bigr) \bigl(
\mathrm{E}\mathbf{A}^{-1}(z) \bigr)\mathbf{x}_n
\\
&=& \bigl(z\underline{m}(z)+z \bigr)^{-1}\mathbf{x}_n^*
\Biggl[\sum_{j=1}^N\mathrm{E}
\beta_j(z) \biggl(\mathbf{r}_j\mathbf{r}_j^*
\mathbf{A}_j^{-1}(z)-\frac{1}{N}\mathrm{E}
\mathbf{A}^{-1}(z) \biggr) \Biggr]\mathbf{x}_n
\\
&&{} +m(z) \bigl(z\mathrm{E}\underline{m}_n(z)-z\underline{m}(z)\bigr)\mathrm{E}m_n^H(z)
\\
&=:&\delta_1+\delta_2,
\end{eqnarray*}
where
\begin{eqnarray*}
\delta_1&=& \bigl(z\underline{m}(z)+z \bigr)^{-1}
\mathbf{x}_n^* \Biggl[\sum_{j=1}^N
\mathrm{E}\beta_j(z) \biggl(\mathbf{r}_j
\mathbf{r}_j^*\mathbf{A}_j^{-1}(z)-
\frac{1}{N}\mathrm{E}\mathbf{A}^{-1}(z) \biggr) \Biggr]\mathbf{x}_n,
\\
\delta_2&=&m_y(z) \bigl(z\mathrm{E}\underline{m}_n(z)-z
\underline{m}(z) \bigr)\mathrm{E}m_n^H(z).
\end{eqnarray*}

%
%
\begin{lem}\label{deltaH}
If
\[
|\delta_1|\leq\frac{C_1}{Nvv_y} \biggl(\frac{1}{v_y}+
\frac{\Delta
^H}{v} \biggr)
\]
holds for some constants $C_0$ and $C_1$, when $v^2v_y\geq C_0N^{-1}$,
under the conditions of Theorem \ref{theorem1}, there exists a constant
$C$ such that $\Delta^H\leq Cv/v_y$.
\end{lem}

\begin{pf}
According to Lemma \ref{main},
\begin{eqnarray*}
\Delta^H&\leq& K_1\int_{-A}^A\bigl|
\mathrm{E}m_n^H(z)-m_y(z)\bigr|\,du
\\
&&{} +K_2v^{-1}\int_{|x|>B}\bigl|\mathrm{E}H^{\mathbf{S}_n}(x)-F_{y_n}(x)\bigr|\,dx
\\
&&{}+K_3v^{-1}\sup_x\int
_{|t|<v}\bigl|F_{y_n}(x+t)-F_{y_n}(x)\bigr|\,dt.
\end{eqnarray*}
From Lemmas \ref{B1} and \ref{Ap7}, we know that there exists a
positive constant $C$, such that
%
%
\begin{eqnarray}\label{2n3}
&&K_2v^{-1}\int_{|x|>B}\bigl|
\mathrm{E}H^{\mathbf{S}_n}(x)-F_{y_n}(x)\bigr|\,dx
\nonumber
\\
&&\qquad\quad{}+K_3v^{-1}\sup_x\int
_{|t|<v}\bigl|F_{y_n}(x+t)-F_{y_n}(x)\bigr|\,dt
\\
&&\qquad \leq Cv/v_y.\nonumber
\end{eqnarray}
Notice that
\begin{eqnarray*}
&& \int_{-A}^A\bigl|\mathrm{E}m_n^H(z)-m_y(z)\bigr|\,du
\\
&&\qquad \leq\int_{-A}^A|\delta_1|\,du+\int_{-A}^A|\delta_2|\,du
\\
&&\qquad \leq\int_{-A}^A|\delta_1|\,du+\int
_{-A}^A\bigl|zm_y(z)\bigr|\bigl|\mathrm{E}m_n^H(z)\bigr|\bigl|\mathrm{E}\underline{m}_n(z)-
\underline{m}(z)\bigr|\,du
\\
&&\qquad \leq\int_{-A}^A|\delta_1|\,du+\int
_{-A}^A\bigl|zm_y(z)\bigr|\bigl|\mathrm{E}m_n^H(z)-m_y(z)\bigr|\bigl|\mathrm{E}
\underline{m}_n(z)-\underline{m}(z)\bigr|\,du
\\
&&\quad\qquad{} +\int_{-A}^A\bigl|zm_y(z)\bigr|\bigl|m_y(z)\bigr|\bigl|
\mathrm{E}\underline{m}_n(z)-\underline{m}(z)\bigr|\,du.
\end{eqnarray*}
Lemma \ref{Ap0}, (\ref{mnmn}) and (\ref{AA}) imply that
\[
\bigl|\mathrm{E}\underline{m}_n(z)-\underline{m}(z)\bigr|=|y_n|\bigl|
\mathrm{E}m_n(z)-m_y(z)\bigr|\leq\frac{C}{Nv^{3/2}v_y^2}
\]
and
\[
\int_{-A}^A\bigl|\mathrm{E}\underline{m}_n(z)-
\underline{m}(z)\bigr|\,du=\int_{-A}^A|y_n|\bigl|
\mathrm{E}m_n(z)-m_y(z)\bigr|\,du\leq Cv.
\]
From (2.3) in \cite{bai1993b}, we have
\begin{eqnarray*}
zm_y(z)&=&\frac{1-y-z+\sqrt{(1+y-z)^2-4y}}{2y}
\\
&=&-1-\frac{1}{2\sqrt{y}}m_{\mathrm{semi}} \biggl(\frac{z-1-y}{\sqrt
{y}} \biggr),
\end{eqnarray*}
where $m_{\mathrm{semi}}(\cdot)$ denotes the Stieltjes transform of the
semicircle law, see (3.2) in~\cite{bai1993a}. Therefore $|zm_y(z)|$ is
bounded by a constant, for $|m_{\mathrm{semi}}(\cdot)|\leq1$, see (3.3)
in~\cite{bai1993a}.

Combined with Lemma \ref{B7}, there exist constants $C_2$, $C_3$, such
that
\begin{eqnarray*}
&&\int_{-A}^A\bigl|\mathrm{E}m_n^H(z)-m_y(z)\bigr|\,du
\\
&&\qquad \leq\int_{-A}^A|\delta_1|\,du +
\frac{C_2}{Nv^{3/2}v_y^2}\int_{-A}^A\bigl|\mathrm{E}m_n^H(z)-m_y(z)\bigr|\,du+\frac{C_3v}{v_y}.
\end{eqnarray*}
Given $v^2v_y\geq C_0N^{-1}$, for $v_y\geq\sqrt{v}$, we have
$v^{3/2}v_y^2\geq v^2v_y\geq C_0N^{-1}$. For a large enough $C_0$ such
that $C_2/C_0\leq1/2$, we have
\[
\int_{-A}^A\bigl|\mathrm{E}m_n^H(z)-m_y(z)\bigr|\,du
\leq\int_{-A}^A|\delta_1|\,du+
\frac{Cv}{v_y}.
\]
As $|\delta_1|\leq\frac{C_1}{Nvv_y} (\frac{1}{v_y}+\frac{\Delta ^H}{v}
)$ and $v^2v_y\geq C_0N^{-1}$, if $C_0\geq4AC_1K_1$, we then have
%
%
\begin{eqnarray}\label{1}
\int_{-A}^A\bigl|\mathrm{E}m_n^H(z)-m_y(z)\bigr|\,du
&\leq&\frac{2AC_1}{Nv^2v_y}\frac{v}{v_y}+\frac{2AC_1}{Nv^2v_y}\Delta
^H+\frac{Cv}{v_y}
\nonumber\\[-8pt]\\[-8pt]
&\leq&\frac{\Delta^H}{2K_1}+\frac{Cv}{v_y}.\nonumber
\end{eqnarray}
Thus, from Lemma \ref{main}, equations (\ref{1}) and (\ref{2n3}), we
conclude that there exists a~constant $C$, such that
\[
\Delta^H\leq\frac{Cv}{v_y}.
\]
The proof is complete.
\end{pf}

To finish the proof of Theorem \ref{theorem1}, we choose $v=\frac
{\sqrt{C_0N^{-1}}}{\sqrt{\sqrt{a}+N^{-1/4}}}$ such that $
v^2v_y=C_0N^{-1}\frac{v_y}{\sqrt{a}+N^{-1/4}}\geq C_0N^{-1}. $
According to Lemma \ref{deltaH}, we know that
\[
\Delta^H\leq\frac{Cv}{v_y}\leq CN^{-1/2} \bigl(
\sqrt{a}+N^{-1/4} \bigr)^{-3/2}.
\]

\noindent If $\sqrt{a}<N^{-1/4}$, $\Delta^H\leq CN^{-1/2} (N^{-1/4}
    )^{-3/2}=O(N^{-1/8})$.

\noindent If $\sqrt{a}\geq N^{-1/4}$, $\Delta^H\leq CN^{-1/2} (\sqrt{a}
    )^{-3/2}=O(N^{-1/2}a^{-3/4})$.

Thus, the proof of Theorem \ref{theorem1} is complete.

\section{\texorpdfstring{The bound for $|\delta_1|$}
{The bound for |delta1|}}\label{sec4}
In this section, we are going to show that when $v^2v_y\geq C_0N^{-1}$,
$|\delta_1|$ is indeed bounded by $\frac{C_1}{Nvv_y} (\frac
{1}{v_y}+\frac{\Delta^H}{v} )$, as required by Lemma \ref{deltaH}.

From $\delta_1$, we can further write $\delta_1=\delta_{11}+\delta
_{12}+\delta_{13}$, where
\begin{eqnarray*}
\delta_{11}&=&N \bigl(z\underline{m}(z)+z \bigr)^{-1}
\mathrm{E} \biggl[\beta_1(z) \biggl(\mathbf{r}_1^*
\mathbf{A}_1^{-1}(z)\mathbf{x}_n \mathbf
{x}_n^* \mathbf{r}_1-\frac{1}{N}
\mathbf{x}_n^* \mathbf{A}_1^{-1}(z)
\mathbf{x}_n \biggr) \biggr]
\\
&=&N \bigl(z\underline{m}(z)+z \bigr)^{-1} \mathrm{E} \bigl(
\beta_1(z)\alpha_1(z) \bigr),
\\
\delta_{12}&=& \bigl(z\underline{m}(z)+z \bigr)^{-1}
\mathrm{E} \bigl[\beta_1(z)\mathbf{x}_n^* \bigl(
\mathbf{A}_1^{-1}(z)-\mathbf{A}^{-1}(z) \bigr)
\mathbf{x}_n \bigr],
\\
\delta_{13}&=& \bigl(z\underline{m}(z)+z \bigr)^{-1}
\mathrm{E} \bigl[\beta_1(z)\mathbf{x}_n^* \bigl(
\mathbf{A}^{-1}(z)-\mathrm{E}\mathbf{A}^{-1}(z) \bigr)
\mathbf{x}_n \bigr].
\end{eqnarray*}

According to (\ref{underlinemm}) and Lemma \ref{B7},
%
%
\begin{equation}\label{mz}
\bigl| \bigl(z\underline{m}(z)+z
\bigr)^{-1} \bigr|= \bigl|-m_y(z) \bigr|\leq \frac{C}{v_y}
\end{equation}
for some constant $C$. Using identity (\ref{betabxi}) three times, we have
\[
\beta_1(z)=b_1(z)-b_1^2(z)
\xi_1(z)+b_1^3(z)\xi_1^2(z)-b_1^3(z)
\beta_1(z)\xi_1^3(z).
\]
Notice that $\mathrm{E}\alpha_1(z)=0$ and $b_1(z)$ is bounded by a
constant (due to Lem\-ma~\ref{Ap2}), we then have
%
%
\begin{eqnarray}\label{eq33}
\qquad|\delta_{11}|&\leq& \frac{CN}{v_y} \bigl|\mathrm{E}\beta_1(z)\alpha_1(z) \bigr|
\nonumber\\[-8pt]\\[-8pt]
&\leq& \frac{CN}{v_y} \bigl( \bigl|\mathrm{E}\xi_1(z)
\alpha_1(z) \bigr|+ \bigl|\mathrm{E}\xi_1^2(z)
\alpha_1(z) \bigr|+ \bigl|\mathrm{E}\beta_1(z)\xi_1^3(z)
\alpha_1(z) \bigr| \bigr).\nonumber
\end{eqnarray}

Let us start with the first term in the above upper bound of $|\delta
_{11}|$ as in (\ref{eq33}). Note that $\mathbf{r}_1$ and $\mathbf
{A}_1^{-1}(z)$ are independent. Therefore, for any integer $p>0$ we
have
\[
\mathrm{E} \bigl(\operatorname{tr}\mathbf{A}_1^{-1}(z)-\mathrm{E}\operatorname{tr}
\mathbf{A}_1^{-1}(z) \bigr)\alpha_1(z)=0
\]
and
\[
\mathrm{E}
\bigl(\operatorname{tr}\mathbf{A}_1^{-1}(z)\bigr)^p\alpha_1(z)=0.
\]
Denote $\mathbf{A}_1^{-1}(z)= (a_{ij} )_{n\times n}$, $\mathbf
{A}_1^{-1}(z)\mathbf{x}_n\mathbf{x}_n^*= (b_{ij} )_{n\times n}$, and
$\mathbf{e}_i$ be the $i$th canonical basis vector, that is, the
$n$-vector whose coordinates are all 0 except that the $i$th coordinate
is 1. Then Lemmas \ref{B2}, \ref{Ap4}, \ref{Ap5} and the inequality $\|
\mathbf{A}_1^{-1}(z)\|\leq1/v$ imply that
\begin{eqnarray*}
&& \bigl|\mathrm{E}\xi_1(z)\alpha_1(z) \bigr|
\\
&&\qquad = \bigl|\mathrm{E}\hat{\xi
}_1(z)\alpha_1(z) \bigr|
\\
&&\qquad \leq\frac{C}{N^2}\mathrm{E} \Biggl(\operatorname{tr} \bigl(\mathbf
{A}_1^{-1}(z)\mathbf{A}_1^{-1}(
\bar{z})\mathbf{x}_n\mathbf{x}_n^* \bigr)+\sum
_{i=1}^na_{ii}b_{ii} \Biggr)
\\
&&\qquad \leq\frac{C}{N^2v} \Biggl\{\frac1{v_y}+\frac{\Delta^H}{v} +v
\Biggl(\sum_{i=1}^N\mathrm{E}\bigl|
\mathbf{x}_n^*\mathbf{A}_1^{-1}(\bar{z})
\mathbf{e}_i\bigr|^2 \sum_{i=1}^N
\mathrm{E}|a_{ii}|^2\bigl|\mathbf{x}_n^*
\mathbf{e}_i\bigr|^2 \Biggr)^{1/2} \Biggr\}
\\
&&\qquad \leq \frac{C}{N^2v} \biggl(\frac1{v_y}+\frac{\Delta^H}{v}\biggr).
\end{eqnarray*}
In the above, we use the following two results, which can be proved by
applying Lemmas \ref{Ap4} and \ref{Ap5}:
\begin{eqnarray*}
\mathrm{E}|a_{11}|&=&\mathrm{E}\bigl|\mathbf{e}_1^*
\mathbf{A}_1^{-1}(z)\mathbf{e}_1\bigr| \leq C
\biggl(\frac1{v_y}+\frac{\Delta^H}{v} \biggr),
\\
\mathrm{E}|a_{11}|^2&\leq&\mathrm{E}\bigl|\mathbf{e}_1^*
\mathbf{A}_1^{-1}(z)\mathbf{A}_1^{-1}(
\bar{z})\mathbf{e}_1\bigr| \leq\frac{C}{v} \biggl(
\frac1{v_y}+\frac{\Delta^H}{v} \biggr).
\end{eqnarray*}
Hence, we have shown that
%
%
\begin{eqnarray}\label{delta111}
\bigl|\mathrm{E}\xi_1(z)\alpha_1(z) \bigr|&\leq&
\frac{C}{N^2v} \biggl(\frac1{v_y}+\frac{\Delta^H}{v} \biggr).
\end{eqnarray}

Let us denote $\mathbf{X}_1^*$ the conjugate transpose of $\mathbf
{X}_1$, that is, $\mathbf{X}_1^*= (\bar{X}_{11},\bar{X}_{21},\break
\ldots,\bar{X}_{n1} )$. Then we can rewrite the second term in the
upper bound of $|\delta_{11}|$ as
\begin{eqnarray*}
\mathrm{E}\xi_1^2(z)\alpha_1(z)&=&
\frac{1}{N^3}\mathrm{E} \biggl[ \biggl(\sum_{i\ne j}a_{ij}
\bar{X}_{i1}X_{j1} +\sum_i
\bigl(a_{ii}\bigl|X_{1i}^2\bigr|-\mathrm{E}a_{ii}
\bigr) \biggr)^2
\\
&&\hspace*{89pt}{}\times\biggl(\sum_{i,j}b_{ij}(
\bar{X}_{i1}X_{j1}-\delta_{ij}) \biggr) \biggr]
\\
&=&\frac{1}{N^3}\mathrm{E} \biggl[ \biggl\{ \biggl(\sum
_{i\ne j}a_{ij}\bar{X}_{i1}X_{j1}
\biggr)^2+ 2 \biggl(\sum_{i}
\bigl(a_{ii}\bigl|X_{1i}^2\bigr|-\mathrm{E}a_{ii}
\bigr) \biggr)
\\
&&\hspace*{33pt}{}\times\biggl(\sum_{i\ne j}a_{ij}
\bar{X}_{i1}X_{j1} \biggr)+ \biggl(\sum
_{i}a_{ii}\bigl|X_{1i}^2\bigr|-
\mathrm{E}a_{ii} \biggr)^2 \biggr\}
\\
&&\hspace*{64pt}{}\times\biggl\{\sum_{i\ne j}b_{ij}
\bar{X}_{i1}X_{j1}+\sum_{i}b_{ii}
\bigl(\bigl|X_{1i}^2\bigr|-1\bigr) \biggr\} \biggr]
\\
&\leq&\frac{C}{N^3}\mathrm{E} \biggl(\sum_i\bigl|a_{ii}^2b_{ii}\bigr|+
\sum_i|\mathrm{E}a_{ii}|^2|b_{ii}|
+\sum_{i\neq j}\bigl|a_{ij}^2b_{ii}\bigr|
\\
&&\hspace*{28pt}{}+ \sum_{i\neq j}\bigl|a_{ij}^2b_{ij}\bigr|+
\sum_{i\neq j}|a_{ij}a_{ii}b_{jj}|
+\sum_{i\neq j}|a_{ii}a_{ij}b_{ij}|
\biggr)
\\
&&{} +\frac{C}{N^3}\sum_{\iota,\tau}\biggl\vert
\sum_{i\neq j\ne k}\mathrm{E}a_{ij}^{\iota}a_{jk}^{\tau
}b_{ik}\biggr\vert,
\end{eqnarray*}
where $a_{ij}^{\iota}$ and $a_{ij}^{\tau}$ denote $a_{ij}$ or $\bar
a_{ij}$. By following the similar proofs in
establishing~(\ref{delta111}), we are able to show that
\begin{eqnarray*}
\mathrm{E}\sum_i\bigl|a_{ii}^2b_{ii}\bigr|&
\leq&\frac1{v}\sum_i\mathrm{E}|a_{ii}b_{ii}|
\\
&\leq& \frac1v \bigl(\mathrm{E}\bigl|a_{11}^2\bigr|\mathrm{E}\mathbf{x}_n^*
\mathbf{A}_1^{-1}(z)\mathbf{A}_1^{-1}(\bar z)
\mathbf{x}_n \bigr)^{1/2}
\\
&\leq& \frac{C}{v^{2}} \biggl(\frac1{v_y}+\frac{\Delta^H}{v}
\biggr),
\\
\mathrm{E}\sum_i|\mathrm{E}a_{ii}|^2|b_{ii}|&
\leq&|\mathrm{E}a_{11}|^2 \bigl(\mathrm{E}\mathbf{x}_n^*
\mathbf{A}_1^{-1}(z)\mathbf{A}_1^{-1}(\bar z)
\mathbf{x}_n \bigr)^{1/2}
\\
&\leq& \biggl[\frac{C}{v} \biggl(\frac1{v_y}+
\frac{\Delta^H}{v} \biggr) \biggr]^{3/2}.
\end{eqnarray*}

The Cauchy--Schwarz inequality implies that
\begin{eqnarray*}
\mathrm{E}\sum_{i\ne j}\bigl|a_{ij}^2b_{ii}\bigr|&
\leq& \mathrm{E} \biggl(\sum_i|b_{ii}|
\biggl(\sum_j\bigl|a_{ij}^2\bigr|
\biggr) \biggr)
\\[-1.8pt]
&\leq& \biggl(\sum_i\mathrm{E}\bigl|
\mathbf{x}_n^*\mathbf{A}_1^{-1}(z)\mathbf
{e}_i\bigr|^2 \biggr)^{1/2} \biggl( \sum
_i\bigl|\mathbf{x}_n^*\mathbf{e}_i\bigr|^2
\mathrm{E} \biggl(\sum_j\bigl|a_{ij}^2\bigr|
\biggr)^2 \biggr)^{1/2}
\\[-1.8pt]
&=& \bigl(\mathrm{E}\mathbf{x}^*_n\bigl(\mathbf{A}_1^{-1}(z)
\mathbf{A}_1^{-1}(\bar{z})\bigr)\mathbf{x}_n
\bigr)^{1/2} \bigl(\mathrm{E}\bigl|\bigl(\mathbf{A}_1^{-1}(z)
\mathbf{A}_1^{-1}(\bar{z})_{11}
\bigr)^2\bigr| \bigr)^{1/2}
\\[-1.8pt]
&\leq&\frac{C}{v^{3/2}} \biggl(\frac1{v_y}+\frac{\Delta^H}{v}
\biggr)^{3/2},
\\[-1.8pt]
\mathrm{E}\sum_{i\ne j}a_{ij}^2b_{ij} &\leq&
\biggl[\sum_{ij}\mathrm{E}\bigl|a_{ij}
\mathbf{x}_n^*\mathbf{A}_1^{-1}(z)\mathbf{e}_i\bigr|^2
\sum_{ij}\mathrm{E}\bigl|a_{ij}\mathbf{x}_n^* \mathbf{e}_j\bigr|^2
\biggr]^{1/2}
\\[-1.8pt]
&\leq& \biggl[\sum_{i}\mathrm{E}\bigl| \bigl(
\mathbf{A}_1^{-1}(z)\mathbf{A}_1^{-1}(\bar{z})
\bigr)_{ii}\bigl(\mathbf{x}_n^*\mathbf{A}_1^{-1}(z)
\mathbf{e}_i\bigr)^2 \bigr|
\\[-1.8pt]
&&\hspace*{21pt}{}\times\sum_{j}\mathrm{E}\bigl|
\bigl(\mathbf{A}_1^{-1}( \bar{z})\mathbf{A}_1^{-1}(z) \bigr)_{jj}\bigl(
\mathbf{x}_n^*\mathbf{e}_j\bigr)^2 \bigr| \biggr]^{1/2}
\\[-1.8pt]
&\leq& \frac{C}{v^2} \biggl(\frac1{v_y}+\frac{\Delta^H}{v}
\biggr)^{3/2},
\\[-1.8pt]
\mathrm{E}\sum_{i\ne j}|a_{ii}a_{ij}b_{jj}|& \leq&
\biggl[\sum_{ij}\mathrm{E}\bigl|a_{ii}
\mathbf{x}_n^*\mathbf{A}_1^{-1}(z)\mathbf{e}_j\bigr|^2
\sum_{ij}\mathrm{E}\bigl|a_{ij}\mathbf{x}_n^* \mathbf{e}_j\bigr|^2
\biggr]^{1/2}
\\[-1.8pt]
&\leq& \biggl[\sum_{i}\mathrm{E}\bigl|a_{ii}^2\bigr|
\mathbf{x}_n^*\mathbf{A}_1^{-1}(z)\mathbf{A}_1^{-1}(
\bar{z})\mathbf{x}_n
\\[-1.8pt]
&&\hspace*{5pt}{}\times\sum_{j} \mathrm{E} \bigl(
\mathbf{A}_1^{-1}(z)\mathbf{A}_1^{-1}(\bar{z})
\bigr)_{jj}\bigl|\mathbf{x}_n^*\mathbf{e}_j\bigr|^2 \biggr]^{1/2}
\\[-1.8pt]
&\leq&C\frac{\sqrt{N}}{v^2}
\biggl(\frac1{v_y}+\frac{\Delta^H}{v}\biggr),
\\[-1.8pt]
\mathrm{E}\sum_{i\neq j}|a_{ii}a_{ij}b_{ij}|& \leq&
\biggl[\sum_{ij}\mathrm{E}\bigl|a_{ii}
\mathbf{x}_n^*\mathbf{A}_1^{-1}\mathbf{e}_j\bigr|^2
\sum_{ij}\mathrm{E}\bigl|a_{ij}\mathbf{x}_n^* \mathbf{e}_i\bigr|^2
\biggr]^{1/2}
\\[-1.8pt]
&\leq&C\frac{\sqrt{N}}{v^2} \biggl(\frac1{v_y}+\frac{\Delta^H}{v}
\biggr).
\end{eqnarray*}
Finally, we establish
\begin{eqnarray*}
&&\sum_{\iota,\tau}\biggl\vert\sum
_{i\neq j\ne k}\mathrm{E}a_{ij}^{\iota
}a_{jk}^{\tau}b_{ik}
\biggr\vert\leq C\frac{\sqrt{N}}{v^2} \biggl(\frac1{v_y}+
\frac{\Delta^H}{v} \biggr).
\end{eqnarray*}
By inclusive--exclusive principle and what we have just proved, it
remains to show that
%
%
\begin{eqnarray}\label{eqb1}
\sum_{\iota,\tau}\biggl\vert\sum
_{i,j, k}\mathrm{E}a_{ij}^{\iota}a_{jk}^{\tau
}b_{ik}
\biggr\vert&\leq& C\frac{\sqrt{N}}{v^2} \biggl(\frac1{v_y}+
\frac{\Delta^H}{v} \biggr).
\end{eqnarray}
Notice that $\gamma_{ik}=\sum_ja_{ij}a_{jk}$ is the $(i,k)$-element of
$\mathbf{A}_1^{-2}(z)$. We obtain
\begin{eqnarray*}
\biggl\vert\sum_{i,j, k}\mathrm{E}a_{ij}a_{jk}b_{ik}
\biggr\vert&=&\biggl\vert\sum_{i,k}\gamma_{ik}b_{ik}\biggr\vert
\\
&\leq& \biggl(\sum_{i,k}\mathrm{E}|\gamma_{ik}|^2
\sum_{i,k}\mathrm{E}|b_{ik}|^2
\biggr)^{1/2}
\\
&=& \bigl(\mathrm{E}\operatorname{tr} \bigl(\mathbf{A}_1^{-2}(z)
\mathbf{A}_1^{-2}(\bar{z}) \bigr)\mathrm{E}\mathbf{x}_n^*
\mathbf{A}_1^{-1}(z)\mathbf{A}_1^{-1}(\bar{z})
\mathbf{x}_n \bigr)^{1/2}
\\
&\leq&C\frac{\sqrt{N}}{v^2} \biggl(\frac1{v_y}+\frac{\Delta^H}{v}
\biggr).
\end{eqnarray*}
Similarly, one can prove the other terms of (\ref{eqb1}) share this
common bound.

In summary, when $v^2v_y\geq C_0N^{-1}$ it holds that
%
%
\begin{equation}\label{delta112}
\bigl|\mathrm{E}\xi_1^2(z)\alpha_1(z) \bigr|\leq \frac{C}{N^2v}
\biggl(\frac1{v_y}+\frac{\Delta^H}{v} \biggr).
\end{equation}

For the last term in the upper bound of $|\delta_{11}|$, we apply Lemma
\ref{Ap3} and the Cauchy--Schwarz inequality again. In particular, for
any fixed $t>0$, we have that
%
%
\begin{eqnarray}\label{delta113}
\bigl|\mathrm{E}\beta_1(z)\xi_1^3(z)
\alpha_1(z) \bigr| &\leq& C\mathrm{E} \bigl|\xi_1^3(z)
\alpha_1(z) \bigr|+o\bigl(N^{-t}\bigr)
\nonumber
\\
&\leq& C \bigl(\mathrm{E} \bigl|\xi_1(z) \bigr|^6
\bigr)^{1/2} \bigl( \mathrm{E} \bigl|\alpha_1(z) \bigr|^2
\bigr)^{1/2}+o\bigl(N^{-t}\bigr)
\\
&\leq&\frac{C}{N^{5/2}v^2v_y^{3/2}} \biggl(\frac1{v_y}+\frac
{\Delta^H}{v}\biggr)^{1/2}.\nonumber
\end{eqnarray}
The last inequality in (\ref{delta113}) is due to Lemmas \ref{Ap1} and
\ref{Ap6}. Therefore, for any $v^2v_y \geq C_0N^{-1}$, (\ref{delta111}),
(\ref{delta112}) and (\ref{delta113}) lead us to
%
%
\begin{equation}\label{delta11}
|\delta_{11}| \leq\frac{C}{Nvv_y}
\biggl(\frac1{v_y}+ \frac{\Delta ^H}{v} \biggr).
\end{equation}

To establish the upper bound for $|\delta_{12}|$, we will make use of
the following equality:
%
%
\begin{equation}\label{aa1}
\mathbf{A}_1^{-1}(z)-\mathbf{A}^{-1}(z)= \beta_1(z)\mathbf{A}_1^{-1}(z)
\mathbf{r}_1 \mathbf{r}_1^*\mathbf{A}_1^{-1}(z).
\end{equation}
Note that (\ref{mz}) implies that
%
%
\begin{eqnarray}\label{delta22}
\qquad\quad|\delta_{12}| &\leq& \frac{C}{v_y} \bigl|\mathrm{E}
\beta_1(z)\mathbf{x}_n^* \bigl(\mathbf{A}_1^{-1}(z)-
\mathbf{A}^{-1}(z) \bigr)\mathbf{x}_n \bigr|
\nonumber
\\
&\leq& \frac{C}{v_y} \bigl|\mathrm{E}\beta_1^2(z)
\mathbf{x}_n^*\mathbf{A}_1^{-1}(z)
\mathbf{r}_1 \mathbf{r}_1^* \mathbf{A}_1^{-1}(z)
\mathbf{x}_n \bigr|
\nonumber
\\
&\leq& \frac{C}{Nv_y} \mathrm{E} \bigl|\mathbf{X}_1^*\mathbf
{A}_1^{-1}(z)\mathbf{x}_n
\mathbf{x}_n^*\mathbf{A}_1^{-1}(z)\mathbf
{X}_1 \bigr|+o\bigl(N^{-t}\bigr)\quad(\mbox{see Lemma~\ref{Ap3}})
\nonumber\\[-8pt]\\[-8pt]
&\leq& \frac{C}{Nv_y} \bigl|\mathrm{E}\operatorname{tr} \bigl(\mathbf{A}_1^{-1}(z)
\mathbf{x}_n \mathbf{x}_n^*\mathbf{A}_1^{-1}(z)
\bigr) \bigr|\qquad(\mbox{see Lemma \ref{B4}})\nonumber
\\
&=& \frac{C}{Nv_y} \bigl|\mathrm{E}\mathbf{x}_n^*\mathbf
{A}_1^{-2}(z)\mathbf{x}_n \bigr|\nonumber
\\
&\leq& \frac{C}{Nvv_y} \biggl(\frac1{v_y}+\frac{\Delta^H}{v}
\biggr).\nonumber
\end{eqnarray}

At last, we establish the upper bound for $\delta_{13}$.

By (\ref{mz}), Lemma \ref{Ap2}, and the fact
%
%
\begin{equation}\label{bbeta}
\beta_1(z)=b_1(z)-b_1(z) \beta_1(z)\xi_1(z),
\end{equation}
we obtain
\begin{eqnarray*}
|\delta_{13}| &\leq& \frac{C}{v_y} \bigl|\mathrm{E} \bigl\{
\beta_1(z)\mathbf{x}_n^* \bigl(\mathbf{A}^{-1}(z)-
\mathrm{E}\mathbf{A}^{-1}(z) \bigr)\mathbf{x}_n \bigr\} \bigr|
\\
&=& \frac{C}{v_y} \bigl|\mathrm{E} \bigl\{b(z)\beta_1(z)
\xi_1(z)\mathbf{x}_n^* \bigl(\mathbf{A}^{-1}(z)-
\mathrm{E}\mathbf{A}^{-1}(z) \bigr)\mathbf{x}_n \bigr\} \bigr|
\\
&\leq& \frac{C}{v_y} \bigl|\mathrm{E} \bigl\{\xi_1(z)
\mathbf{x}_n^* \bigl(\mathbf{A}^{-1}(z)-\mathrm{E}
\mathbf{A}^{-1}(z) \bigr)\mathbf{x}_n \bigr\} \bigr|
\\
&&{} +\frac{C}{v_y} \bigl|\mathrm{E} \bigl\{\beta_1(z)
\xi_1^2(z)\mathbf{x}_n^* \bigl(
\mathbf{A}^{-1}(z)-\mathrm{E}\mathbf{A}^{-1}(z) \bigr)\mathbf
{x}_n \bigr\} \bigr|.
\end{eqnarray*}
From $C_r$ inequality (see Lo{\`e}ve \citet{Lo78}), we have
$|\delta_{13}|\leq\frac{C}{v_y}
(|\mathrm{II}_1|+|\mathrm{II}_2|+|\mathrm{II}_3| )$, where
\begin{eqnarray*}
\mathrm{II}_1&=&\mathrm{E} \bigl\{\xi_1(z)
\mathbf{x}_n^* \bigl(\mathbf{A}_1^{-1}(z)-
\mathrm{E}\mathbf{A}_1^{-1}(z) \bigr)\mathbf{x}_n
\bigr\},
\\
\mathrm{II}_2&=&\mathrm{E} \bigl\{\xi_1(z)
\mathbf{x}_n^* \bigl(\mathbf{A}^{-1}(z)-
\mathbf{A}_1^{-1}(z) \bigr)\mathbf{x}_n \bigr
\},
\\
\mathrm{II}_3&=& \mathrm{E} \bigl\{\beta_1(z)
\xi_1(z)\mathbf{x}_n^*\mathrm{E} \bigl(
\mathbf{A}^{-1}(z)-\mathbf{A}_1^{-1}(z) \bigr)
\mathbf{x}_n \bigr\}.
\end{eqnarray*}

It should be noted that $\mathrm{E} (\xi_1(z) \mid\mathbf
{A}_1^{-1}(z) )=0$, and $\mathbf{r}_1$ and $\mathbf{A}_1^{-1}(z)$ are
independent. Then we have
\begin{eqnarray*}
\mathrm{II}_1&=&\mathrm{E} \bigl[\mathrm{E} \bigl\{ \xi_1(z)
\mathbf{x}_n^*\bigl(\mathbf{A}_1^{-1}(z)-
\mathrm{E}\mathbf{A}_1^{-1}(z)\bigr)\mathbf{x}_n
\mid\mathbf{A}_1^{-1}(z) \bigr\} \bigr]= 0.
\end{eqnarray*}
By the results in (\ref{aa1}) and (\ref{bbeta}), we have
\begin{eqnarray*}
|\mathrm{II}_2|&=& \bigl|\mathrm{E}\xi_1(z)\beta_1(z)
\mathbf{x}_n^*\mathbf{A}_1^{-1}(z)
\mathbf{r}_1\mathbf{r}_1^*\mathbf{A}_1^{-1}(z)
\mathbf{x}_n \bigr|
\\
&\leq& \bigl|b_1(z) \bigr| \bigl|\mathrm{E}\xi_1(z)
\mathbf{r}_1^*\mathbf{A}_1^{-1}(z)
\mathbf{x}_n\mathbf{x}_n^*\mathbf{A}_1^{-1}(z)
\mathbf{r}_1 \bigr|
\\
&&{}+ \bigl|b_1(z) \bigr| \bigl|\mathrm{E}\beta_1(z)\xi_1^2(z)
\mathbf{r}_1^*\mathbf{A}_1^{-1}(z)
\mathbf{x}_n\mathbf{x}_n^*\mathbf{A}_1^{-1}(z)
\mathbf{r}_1 \bigr|
\\
&=:& \mathrm{III}_1+\mathrm{III}_2,
\end{eqnarray*}
where
\begin{eqnarray*}
\mathrm{III}_1&=& \bigl|b_1(z) \bigr| \biggl|\mathrm{E}
\xi_1(z) \biggl( \mathbf{r}_1^*\mathbf{A}_1^{-1}(z)
\mathbf{x}_n\mathbf{x}_n^*\mathbf{A}_1^{-1}(z)
\mathbf{r}_1- \frac{1}{N}\mathbf{x}_n^*
\mathbf{A}_1^{-2}(z)\mathbf{x}_n \biggr) \biggr|
\\
&\leq& \frac{C}{N^2} \biggl| \mathrm{E}\operatorname{tr} \bigl(\mathbf{x}_n^*\mathbf
{A}_1^{-3}(z)\mathbf{x}_n \bigr)+\mathrm{E}\sum
_{i} \bigl( \mathbf{A}_1^{-1}(z)
\mathbf{x}_n\mathbf{x}_n^*\mathbf{A}_1^{-1}(z)
\bigr)_{ii} a_{ii} \biggr|
\\
&=& \frac{C}{N^2v^2} \biggl(\frac1{v_y}+\frac{\Delta^H}{v}
\biggr).
\end{eqnarray*}
The above inequality follows from Lemmas \ref{B2} and \ref{Ap2}. By
Lemma \ref{Ap3} and the Cauchy--Schwarz inequality, it holds that
\begin{eqnarray*}
\mathrm{III}_2&\leq& C\mathrm{E} \bigl|\xi_1^2(z)
\mathbf{r}_1^*\mathbf{A}_1^{-1}(z)
\mathbf{x}_n\mathbf{x}_n^*\mathbf{A}_1^{-1}(z)
\mathbf{r}_1 \bigr|
\\
&\leq& C \bigl(\mathrm{E}\bigl|\xi_1(z)\bigr|^4
\bigr)^{1/2} \bigl(\mathrm{E}\bigl|\mathbf{r}_1^*
\mathbf{A}_1^{-1}(z)\mathbf{x}_n
\mathbf{x}_n^*\mathbf{A}_1^{-1}(z)
\mathbf{r}_1\bigr|^2 \bigr)^{1/2}
\\
&\leq& C \biggl(\frac{1}{N^2v^2v_y^2}\frac{1}{N^2}\mathrm{E}\bigl|\mathbf
{x}_n^*A_1^{-2}(z)\mathbf{x}_n\bigr|^2
\biggr)^{1/2}\qquad(\mbox{see Lemmas \ref{Ap1} and \ref{B4}})
\\
&\leq& \frac{C}{N^2v^2v_y} \biggl(\frac1{v_y}+\frac{\Delta^H}{v}\biggr).
\end{eqnarray*}
Hence, we have shown that
\[
|\mathrm{II}_2| \leq\frac{C}{N^2v^2v_y} \biggl(\frac1{v_y}+
\frac{\Delta
^H}{v} \biggr).
\]

Moreover, Lemmas \ref{Ap1} and \ref{Ap4}, and the Cauchy--Schwarz
inequality lead us to the following:
\begin{eqnarray*}
|\mathrm{II}_3| &\leq& C \bigl(\mathrm{E} \bigl|\xi_1(z)
\bigr|^2\mathrm{E} \bigl|\mathbf{x}_n^* \bigl(\mathbf{A}^{-1}(z)-
\mathbf{A}_1^{-1}(z) \bigr)\mathbf{x}_n
\bigr|^2 \bigr)^{1/2}
\\
&\leq& \frac{C}{N^{1/2}v^{1/2}v_y^{1/2}} \frac{C}{Nv} \biggl(\frac1{v_y}+
\frac{\Delta^H}{v} \biggr)
\\
&=&\frac{C}{N^{3/2}v^{3/2}v_y^{1/2}} \biggl(\frac1{v_y}+\frac
{\Delta
^H}{v}
\biggr).
\end{eqnarray*}

Therefore, it follows that
%
%
\begin{equation}\label{delta33}
|\delta_{13}| \leq \frac{C}{v_y}
\bigl(|\mathrm{II}_1|+| \mathrm{II}_2|+|\mathrm{II}_3| \bigr)
\leq \frac{C}{N^{3/2}v^{3/2}v_y} \biggl(\frac1{v_y}+\frac{\Delta^H}{v}\biggr).
\end{equation}

As it has been shown in (\ref{delta11}), (\ref{delta22}) and (\ref
{delta33}), we conclude that
%
%
\begin{equation}\label{delta1}
|\delta_1|\leq\frac{C}{Nvv_y} \biggl(\frac1{v_y}+
\frac{\Delta^H}{v} \biggr).
\end{equation}

\section{\texorpdfstring{Proof of Theorem \protect\ref{theorem2}}{Proof of Theorem 1.6}}\label{sec5}
From Lemma \ref{main}, and by replacing $\mathrm{E}H^{\mathbf{S}_n}(x)$ and
$\mathrm{E}m_n^H(z)$ by $H^{\mathbf{S}_n}(x)$ and $m_n^H(z)$, respectively, we
have
\begin{eqnarray*}
\mathrm{E}\Delta_p^H&=:&\mathrm{E}\bigl\|H^{\mathbf{S}_n}(x)-F_{y_n}(x)\bigr\|
\\
&\leq&K_1\int_{-A}^A\mathrm{E}\bigl|m_n^H(z)-\mathrm{E}m_n^H(z)\bigr|\,du
+K_1\int_{-A}^A\bigl|\mathrm{E}m_n^H(z)-m_y(z)\bigr|\,du
\\
&&{}+K_2v^{-1}\int_{|x|>B}\bigl|\mathrm{E}H^{\mathbf{S}_n}(x)-F_{y_n}(x)\bigr|\,dx
\\
&&{}+K_3v^{-1}\sup_x\int_{|t|<v}\bigl|F_{y_n}(x+t)-F_{y_n}(x)\bigr|\,dt
\\
&\leq&K_1\int_{-A}^A\mathrm{E}\bigl|m_n^H(z)-\mathrm{E}m_n^H(z)\bigr|\,du+\Delta^H.
\end{eqnarray*}
As the convergence rate of $\Delta^H$ has already been established in
Theorem \ref{theorem1}, we only focus on the convergence rate of
$\mathrm{E}|m_n^H(z)-\mathrm{E}m_n^H(z)|$.

By Lemma \ref{Ap5} and the Cauchy--Schwarz inequality, it follows that
\begin{eqnarray*}
\mathrm{E} \bigl|m_n^H(z)-\mathrm{E}m_n^H(z)
\bigr| &\leq& \bigl(E\bigl|m_n^H(z)-\mathrm{E}m_n^H(z)\bigr|^2
\bigr)^{1/2}
\\
&\leq& \frac{C}{\sqrt{N}v} \biggl(\frac1{v_y}+\frac{\Delta^H}{v}
\biggr).
\end{eqnarray*}
Together with Lemma \ref{deltaH} that $\Delta^H\leq Cv/v_y$, when
$v^2v_y\geq O(N^{-1})$, we have
\[
\mathrm{E}\bigl\|H^{\mathbf{S}_n}(x)-F_{y_n}(x)\bigr\|\leq\biggl(
\frac{1}{\sqrt{N}v}+v \biggr)\frac{C}{v_y}.
\]
By choosing $v=O(N^{-1/4})$, we obtain
\[
\mathrm{E}\bigl\|H^{\mathbf{S}_n}(x)-F_{y_n}(x)\bigr\|\leq\cases{ O
\bigl(N^{-1/4}a^{-1/2}\bigr), &\quad when $a\geq
N^{-1/4}$, \vspace*{2pt}
\cr
O\bigl(N^{-1/8}\bigr), &\quad
otherwise.}
\]
The proof of Theorem \ref{theorem2} is complete.

\section{\texorpdfstring{Proof of Theorem \protect\ref{theorem3}}{Proof of Theorem 1.8}}\label{sec6}

Notice that the proof of Theorem \ref{theorem3} is almost the same as
that of Theorem \ref{theorem2}.

By Lemma \ref{main}, choosing $v=O(N^{-1/4})$,
\[
\bigl\|H^{\mathbf{S}_n}-F_{y_n} \bigr\| \leq \int_{-A}^A
\bigl|m_n^H(z)-\mathrm{E}m_n^H(z)\bigr|\,du+Cv/v_y.
\]

By Lemma \ref{Ap5}, we have
\[
\mathrm{E}\bigl|m_n^H(z)-\mathrm{E}m_n^H(z)\bigr|^{2l}
\leq CN^{-l}v^{-2l} \biggl(\frac
{1}{v_y}+
\frac{\Delta^H}{v} \biggr)^{2l}.
\]

When $a<N^{-1/4}$, with $v=O(N^{-1/4})$,
\[
N^{2l(1/8-\eta)}\mathrm{E}\bigl|m_n^H(z)-\mathrm{E}m_n^H(z)\bigr|^{2l}
\leq CN^{-2l\eta},
\]
which implies that if we choose an $l$ such that $2l\eta>1$,
\[
\int_{-A}^A \bigl|m_n^H(z)-\mathrm{E}m_n^H(z)\bigr|\,du=o_{\mathrm{a.s.}}\bigl(N^{-1/8+\eta}
\bigr).
\]

When $a\geq N^{-1/4}$, in this case,
by choosing $v=O(N^{-1/4})$, we have
\[
a^l N^{2l(1/4-\eta)}\mathrm{E}\bigl|m_n^H(z)-
\mathrm{E}m_n^H(z)\bigr|^{2l}\leq CN^{-2l\eta}.
\]
Theorem \ref{theorem3} then follows by setting $l>\frac{1}{2\eta}$.
This completes the proof of Theorem~\ref{theorem3}.

\begin{appendix}
\section{}\label{sec7}
In this section, we establish some lemmas which are used in the proofs
of the main theorems.

%
\begin{lem}\label{Ap0}
Under the conditions of Theorem 1.6, for all $|z|<A$ and $v^2v_y\geq
C_0N^{-1}$, we have
\[
\bigl|\mathrm{E}m_n(z)-m_y(z)\bigr|\leq\frac{C}{Nv^{3/2}v_y^2},
\]
where $C_0$ is a constant and $v_y=1-\sqrt{y_n}+\sqrt{v}$.
\end{lem}

\begin{pf}
Since
%
%
\begin{eqnarray}\label{8.3.7}
\mathrm{E}m_n(z)&=&\frac{1}{n}\mathrm{E}\operatorname{tr}(
\mathbf{S}_n-z\mathbf{I}_n)^{-1}
\nonumber
\\
&=&\frac{1}{n}\sum_{k=1}^n
\mathrm{E}\frac{1}{s_{kk}-z-N^{-2}\bolds
{\alpha}_k^*(\mathbf{S}_{nk}-z\mathbf{I}_{n-1})^{-1}\bolds{\alpha
}_k}
\nonumber\\[-8pt]\\[-8pt]
&=&\frac{1}{n}\sum_{k=1}^n
\mathrm{E}\frac{1}{\epsilon
+1-y_n-z-y_nz\mathrm{E}m_n(z)}
\nonumber
\\
&=&-\frac{1}{z+y_n-1+y_nz\mathrm{E}m_n(z)}+\delta_n,\nonumber
\end{eqnarray}
where
\begin{eqnarray*}
s_{kk}&=&\frac{1}{N}\sum_{j=1}^N|X_{kj}|^2,
\\
\mathbf{S}_{nk}&=&\frac{1}{N}\mathbf{X}_{(k)}
\mathbf{X}_{(k)}^*,
\\
\bolds{\alpha}_k&=&\mathbf{X}_{(k)}\mathbf{
\bar{X}}_k,
\\
\epsilon_k&=&(s_{kk}-1)+y_n+y_nz
\mathrm{E}m_n(z)-\frac{1}{N^2}\bolds{\alpha}_k^*(
\mathbf{S}_{nk}-z\mathbf{I}_{n-1})^{-1}\bolds{
\alpha}_k,
\\
\delta_n&=&-\frac{1}{n}\sum_{k=1}^n
b_n\mathrm{E}\beta_k\epsilon_k,
\\
b_n&=&b_n(z)=\frac{1}{z+y_n-1+y_nz\mathrm{E}m_n(z)},
\\
\beta_k&=&\beta_k(z)=\frac{1}{z+y_n-1+y_nz\mathrm{E}m_n(z)-\epsilon_k},
\end{eqnarray*}
where $\mathbf{X}_{(k)}$ is the $(n-1)\times N$ matrix obtained from
$\mathbf{X}$ with its $k$th row removed and $\mathbf{X}_k^*$ is the
$k$th row of $\mathbf{X}$. It has proved that one of the roots of
equation (\ref{8.3.7}) is (see (3.1.7) in \cite{bai1993b})
\[
\mathrm{E}m_n(z)=-\frac{1}{2y_nz} \bigl(z+y_n-1-y_nz
\delta_n- \sqrt{(z+y_n-1+y_nz
\delta_n)^2-4y_nz} \bigr).
\]
The Stieltjes transform of the Mar\v{c}enko--Pastur distribution~with
index $y$ is given by (see
(2.3) in \cite{bai1993b})
\[
m_y(z)=-\frac{y_n+z-1-\sqrt{(1+y_n-z)^2-4y_n}}{2y_nz}.
\]
Thus,
\begin{eqnarray*}
&& \bigl|\mathrm{E}m_n(z)-m_y(z) \bigr|
\\
&&\qquad \leq\frac{|\delta_n|}{2}
\biggl[1+ \frac
{|2(z+y_n-1)-y_nz\delta_n|}{|\sqrt{(z+y_n-1)^2-4y_nz}+\sqrt
{(z+y_n-1+y_nz\delta_n)^2-4y_nz}|} \biggr].
\end{eqnarray*}
Let us define by convention
\[
\Re(\sqrt{z})=\frac{\Im(z)}{\sqrt{2(|z|-\Re(z))}},\qquad
\Im(\sqrt{z})=\frac{|\Im(z)|}{\sqrt{2(|z|+\Re(z))}}.
\]
If $|u-y_n-1| \geq\frac{1}{5(A+1)}$, then the real parts of $\sqrt
{(z+y_n-1)^2-4y_nz}$ and $\sqrt{(z+y_n-1+y_nz\delta_n)^2-4y_nz}$ have
the same sign. Since they both have positive imaginary parts, it
follows that
\begin{eqnarray*}
&& \bigl|\sqrt{(z+y_n-1)^2-4y_nz}+
\sqrt{(z+y_n-1+y_nz\delta_n)^2-4y_nz}
\bigr|
\\
&&\qquad \geq \sqrt{\bigl|\Im\bigl((z+y_n-1)^2-4y_nz
\bigr)\bigr|}=\sqrt{2v(u-y_n-1)} \geq\biggl(\frac{2v}{5(A+1)}
\biggr)^{1/2}.
\end{eqnarray*}
Thus,
%
%
\begin{equation}\label{geq}
\bigl|\mathrm{E}m_n(z)-m_y(z)\bigr|\leq\frac{|\delta_n|}{2} \biggl(1+
\frac{C}{\sqrt{v}} \biggr)\leq\frac{C|\delta_n|}{\sqrt{v}}.
\end{equation}
If $|u-y_n-1|<\frac{1}{5(A+1)}$, we have $|\mathrm{E}m_n(z)-m_y(z)|\leq
C|\delta_n|$.

In \cite{baijack2010} [see the inequality above (8.3.16)], we have
\[
|\delta_n|\leq\frac{C}{Nv^3} \biggl(\Delta+\frac{v}{v_y}
\biggr)^2\leq\frac{C}{Nvv_y^2}.
\]
Combined with (\ref{geq}), we get
\[
\bigl|\mathrm{E}m_n(z)-m_y(z)\bigr|\leq\frac{C}{Nv^{3/2}v_y^2}.
\]
The proof of the lemma is complete.\eject

In addition, the following relevant result which is involved in the
proof of Lemma \ref{deltaH} is presented here:
%
%
\begin{eqnarray} \label{AA}
&& \int_{-A}^A \bigl\vert\mathrm{E}m_n(z)-m_y(z)
\bigr\vert\,du\nonumber
\\
&&\qquad = \int_{|u-y_n-1|\geq 1/(5(A+1)),|u|\leq
A}\bigl\vert\mathrm{E}m_n(z)-m_y(z)\bigr\vert\,du
\nonumber\\[-8pt]\\[-8pt]
&&\quad\qquad{}+\int_{|u-y_n-1|<1/(5(A+1)),|u|\leq A}\bigl\vert\mathrm{E}m_n(z)-m_y(z)
\bigr\vert\,du
\nonumber
\\
&&\qquad \leq Cv.\nonumber
\end{eqnarray}\upqed
\end{pf}

%
%
\begin{lem}\label{Ap1}
Under the conditions of Theorem \ref{theorem2}, for $v^2v_y\geq
C_0N^{-1}$ and $1\leq l \leq3$, there exists a constant $C$, such
that
\[
\mathrm{E} \bigl|\xi_1(z) \bigr|^{2l} \leq\frac{C}{N^lv^lv_y^l}.
\]
\end{lem}

\begin{pf}
By the $C_r$-inequality (see Lo{\`e}ve \citet{Lo78}), it follows that
\begin{eqnarray*}
\mathrm{E} \bigl|\xi_1(z) \bigr|^{2l} &\leq& C \biggl(\mathrm{E} \biggl|
\frac{1}{N}\operatorname{tr}\mathbf{A}_1^{-1}(z)-
\frac
{1}{N}\mathrm{E}\operatorname{tr}\mathbf{A}_1^{-1}(z)
\biggr|^{2l}+ \mathrm{E} \bigl|\hat{\xi}_1(z) \bigr|^{2l} \biggr)
\\
&=:&I_1+I_2.
\end{eqnarray*}
From Lemmas \ref{B6}, \ref{B8}, \ref{B12} and the $C_r$-inequality with
$v^2v_y\geq C_0N^{-1}$, we have
\begin{eqnarray*}\label{xil1}
I_1 &\leq& \mathrm{E}\biggl\vert\frac{1}{N}\operatorname{tr}
\mathbf{A}_1^{-1}(z)-\frac{1}{N}\operatorname{tr}
\mathbf{A}^{-1}(z)\biggr\vert^{2l} +\mathrm{E}\biggl\vert
\frac{1}{N}\operatorname{tr}\mathbf{A}^{-1}(z)-\frac{1}{N}\mathrm{E}\operatorname{tr}\mathbf{A}^{-1}(z)\biggr\vert^{2l}
\\
&&{} + \biggl\vert\frac{1}{N}\mathrm{E}\operatorname{tr}\mathbf{A}^{-1}(z)-
\frac{1}{N}\mathrm{E}\operatorname{tr}\mathbf{A}_1^{-1}(z)\biggr\vert^{2l}
\\
&\leq& C \biggl\{ \biggl(\frac{1}{Nv} \biggr)^{2l}+
\frac
{1}{N^{2l}v^{4l}} \biggl(\Delta+\frac{v}{v_y} \biggr)^l+
\biggl(\frac
{1}{Nv} \biggr)^{2l} \biggr\}
\\
&\leq& \frac{C}{N^{2l}v^{3l}v_y^l}.
\end{eqnarray*}
Under finite 8th moment assumption, for $l\geq2$ we have
\[
\mathrm{E}|X_{11}|^{4l}=\mathrm{E} \bigl\{|X_{11}|^8|X_{11}|^{4l-8}
\mathrm{I}\bigl(|X_{11}|\leq\eta_NN^{1/4}\bigr)
\bigr\}\leq CN^{l-2}.
\]
It can be shown that $\mathbf{B}_1-z\mathbf{I}_n=\mathbf{A}_1$, $ \|
(\mathbf{B}_1-z\mathbf{I}_n )^{-1} \|\leq1/v$, and
\[
\operatorname{tr} \bigl( (\mathbf{B}_1-z\mathbf{I}_n
)^{-1} (\mathbf{B}_1-\bar{z}\mathbf{I}_n
)^{-1} \bigr)=v^{-1}\Im\bigl(\operatorname{tr}(\mathbf{B}_1-z
\mathbf{I}_n)^{-1} \bigr).
\]
Hence, by Lemma \ref{B4},
\begin{eqnarray*}
I_2&=&\frac{1}{N^{2l}}\mathrm{E} \bigl|\mathbf{X}_1^*\mathbf
{A}_1^{-1}(z)\mathbf{X}_1-\operatorname{tr}
\mathbf{A}_1^{-1}(z) \bigr|^{2l}
\\
&\leq& \frac{C}{N^{2l}}\mathrm{E} \bigl\{CN^{l-2}\operatorname{tr} \bigl(
\mathbf{A}_1^{-1}(z)\mathbf{A}_1^{-1}(
\bar{z}) \bigr)^l+ \bigl(C\operatorname{tr}\bigl(\mathbf{A}_1^{-1}(z)
\mathbf{A}_1^{-1}(\bar{z})\bigr) \bigr)^l
\bigr\}
\\
&=& \frac{C}{N^{2l}}\mathrm{E} \bigl\{CN^{l-2}v^{-2l+1}\Im
\bigl(\operatorname{tr}(\mathbf{B}_1-z\mathbf{I}_n)^{-1}
\bigr)
\\
&&\hspace*{40pt}{}+C^lv^{-l} \bigl(\Im\bigl(\operatorname{tr}(
\mathbf{B}_1-z\mathbf{I}_n)^{-1} \bigr)
\bigr)^l \bigr\}
\\
&=& \frac{C}{N^{2l}} \bigl\{CN^{l-1}v^{-2l+1}\mathrm{E} \bigl(
\Im\bigl(m_{F^{\mathbf{B}_1}}(z) \bigr) \bigr)
\\
&&\qquad\qquad{} +C^lN^lv^{-l}\mathrm{E} \bigl(
\Im\bigl(m_{F^{\mathbf{B}_1}}(z) \bigr) \bigr)^l \bigr\}
\\
&\leq& \frac{C}{N^{l+1}v^{2l-1}v_y}+\frac{C}{N^lv^lv_y^l}.
\end{eqnarray*}
The last inequality is due to Lemmas \ref{B6}, \ref{B7}, \ref{B8} and
\begin{eqnarray*}
\bigl|\mathrm{E} \bigl(\Im\bigl(m_{F^{\mathbf{B}_1}}(z)\bigr) \bigr)^l
\bigr| &\leq& \mathrm{E} \bigl|m_{F^{\mathbf{B}_1}}(z)-m_n(z) \bigr|^l +\mathrm{E}
\bigl|m_n(z)-\mathrm{E}m_n(z) \bigr|^l
\\
&&{} + \bigl|\mathrm{E}m_n(z)-m_y(z) \bigr|^l+\bigl|m_y(z) \bigr|^l
\\
&\leq& Cv_y^{-l}.
\end{eqnarray*}
Here let $\Delta= \|EF^{\mathbf{S}_n}-F_{y_n}\|$, by
integration by parts and Lemma \ref{B12}, we have
%
%
\begin{equation}\label{eq611}
\bigl|\mathrm{E}m_n(z)-m_y(z) \bigr|
\leq\frac{C\Delta}{v} \leq \frac{C}{v_y}.
\end{equation}

Therefore, for $1\leq l\leq3$ and $v^2v_y\geq C_0N^{-1}$, it follows
that
\[
\mathrm{E} \bigl|\xi_1(z) \bigr|^l \leq I_1+I_2
\leq\frac{C}{N^lv^lv_y^l}.
\]
This completes the proof.
\end{pf}

%
%
\begin{lem}\label{Ap2}
Under the conditions of Theorem \ref{theorem2}, for all $\bigl|z\big|<A$
and $v^2v_y\geq C_0N^{-1}$, we have
\[
\bigl|b_1(z) \bigr|\leq C.
\]
\end{lem}

\begin{pf}
From (2.3) in \cite{bai1993b}, we have
\[
m_y(z)=\frac{1-y_n-z+\sqrt{(1-y_n-z)^2-4y_nz}}{2y_nz},
\]
where the square root of a complex number is defined as the one with
positive imaginary part.
It then can be verified that
\begin{eqnarray*}
b_0(z)&=:&\frac1{1+y_nm_y(z)} =1+\frac12
\bigl(z-y_n-1-\sqrt{(z-y_n-1)^2-4y_n}
\bigr)
\\
&=&1+\sqrt{y_n}m_{\mathrm{semi}} \biggl(\frac{z-y_n-1}{\sqrt{y_n}}
\biggr),
\end{eqnarray*}
where $m_{\mathrm{semi}}$ denotes the Stieltjes transform of the
semicircular law. As $|m_{\mathrm{semi}}|\leq1$, we conclude that
%
%
\begin{equation}\label{l621}
\bigl\vert b_0(z)\bigr\vert\leq1+\sqrt{y_n}.
\end{equation}
By the relationship between $b_0$ and $b_1$, we have
\[
b_1(z)=\frac{b_0(z)}{1+y_n b_0(z)((1/n) \mathrm{E}\operatorname{tr}\mathbf
{A}^{-1}_1(z)-m_y(z))}.
\]
When $C_0$ is chosen large enough, by Lemma \ref{Ap0}, for all large
$N$ we have
\begin{eqnarray*}
\biggl\vert\frac1n \mathrm{E}\operatorname{tr}\mathbf{A}^{-1}_1(z)-m_y(z)
\biggr\vert&\leq& \frac1n\bigl\vert\mathrm{E}\operatorname{tr} \bigl(\mathbf{A}^{-1}_1(z)-
\mathbf{A}^{-1}(z) \bigr)\bigr\vert+\bigl\vert\mathrm{E}m_n(z)-m_y(z)\bigr\vert
\\
&\leq& \frac1{nv}+\frac{1}{3(1+\sqrt{y_n})}\leq\frac
{2}{3(1+\sqrt{y_n})}
\end{eqnarray*}
and consequently we obtain
\[
\bigl|b_1(z)\bigr|\leq3(1+\sqrt{y_n})\leq C.
\]
Thus, the proof is complete.
\end{pf}

%
\begin{lem}\label{Ap3}
If $ |b_1(z) |\leq C$, then for any fixed $t>0$,
\[
\mathrm{P} \bigl( \bigl|\beta_1(z) \bigr|>2C \bigr)=o\bigl(N^{-t}\bigr).
\]
\end{lem}

\begin{pf}
Note that if $\vert b_1(z)\xi_1(z)\vert\leq1/2$, by Lemma
\ref{Ap2}, we get
\[
\bigl|\beta_1(z) \bigr|=\frac{ |b_1(z) |}{ |1+b_1(z)\xi_1(z) |} \leq\frac{
|b_1(z) |}{1- |b_1(z)
\xi_1(z) |}\leq2
\bigl|b_1(z) \bigr| \leq2C.
\]
As a result,
\begin{eqnarray*}
\mathrm{P} \bigl( \bigl|\beta_1(z) \bigr|>2C \bigr)&\leq& \mathrm{P} \biggl(
\bigl|b_1(z)\xi_1(z) \bigr|>\frac{1}{2} \biggr)
\\
&\leq& \mathrm{P} \biggl( \bigl|\xi_1(z) \bigr|>\frac{1}{2C} \biggr)\qquad (\mbox{see Lemma \ref{Ap2}})
\\
&\leq& (2C )^p\mathrm{E} \bigl|\xi_1(z) \bigr|^p.
\end{eqnarray*}
By the $C_r$-inequality, Lemmas \ref{B8} and \ref{B10}, for some $\eta
=\eta_NN^{-1/4}$ and $p\geq\log N$, we have
\begin{eqnarray*}
\mathrm{E} \bigl|\xi_1(z) \bigr|^p&=& \mathrm{E} \bigl|\hat{\xi}_1(z) \bigr|^p +\mathrm{E}\biggl
\vert\frac{1}{N}\operatorname{tr}\mathbf{A}_1^{-1}(z)-\frac{1}{N}\mathrm{E}\operatorname{tr}
\mathbf{A}_1^{-1}(z)\biggr\vert^p
\\
&\leq& C \bigl(N\eta_N^4N^{-1}
\bigr)^{-1} \bigl(v^{-1}\eta_N^2N^{-1/2}
\bigr)^p+\frac
{C}{N^pv^{3p/2}v_y^{p/2}}
\\
&\leq& C\eta_N^{2p-4}\leq C\eta_N^p.
\end{eqnarray*}
For any fixed $t>0$, when $N$ is large enough so that $\log\eta_N^{-1}>t+1$,
it can be shown that
\begin{eqnarray*}
\mathrm{E}\bigl|\xi_1(z)\bigr|^p &\leq&C e^{-p\log\eta_N^{-1}}
\\
&\leq& C e^{-p(t+1)}
\\
&\leq& Ce^{-(t+1)\log N}
\\
&=&CN^{-t-1}=o\bigl(N^{-t}\bigr).
\end{eqnarray*}
We finish the proof.
\end{pf}

%
%
\begin{lem}\label{Ap4}
If $v^2v_y\geq C_0N^{-1}$, $C_0$ is a large constant. For $l\geq1$, it
holds that
\[
\mathrm{E}\bigl\vert m_{H^{\mathbf{B}_1}}(z)\bigr\vert^{2l} \leq C
\mathrm{E}\bigl\vert m_n^H(z) \bigr\vert
^{2l}.
\]
\end{lem}

\begin{pf}
Recall that
\[
\mathbf{A}_j^{-1}(z)-\mathbf{A}^{-1}(z)=
\beta_j(z)\mathbf{A}_j^{-1}(z)
\mathbf{r}_j\mathbf{r}_j^*\mathbf{A}_j^{-1}(z).
\]
By Lemmas \ref{Ap3} and \ref{B4}, it holds that
\begin{eqnarray*}
&&\mathrm{E}\bigl\vert m_{H^{\mathbf{B}_1}}(z)-m_n^H(z)
\bigr\vert^{2l}
\\
&&\qquad =\mathrm{E}\bigl\vert\mathbf{x}_n^*
\bigl(\mathbf{A}_1^{-1}(z)-\mathbf{A}^{-1}(z)
\bigr)\mathbf{x}_n\bigr\vert^{2l}
\\
&&\qquad =\mathrm{E}\bigl\vert\mathbf{x}_n^*\beta_1(z)
\mathbf{A}_1^{-1}(z)\mathbf{r}_1
\mathbf{r}_1^*\mathbf{A}_1^{-1}(z)
\mathbf{x}_n\bigr\vert^{2l}
\\
&&\qquad =\mathrm{E}\bigl\vert\mathbf{x}_n^*\beta_1(z)
\mathbf{A}_1^{-1}(z)\mathbf{r}_1
\mathbf{r}_1^*\mathbf{A}_1^{-1}(z)
\mathbf{x}_n\bigr\vert^{2l} \mathrm{I}\bigl(\bigl|
\beta_1(z)\bigr|\leq C\bigr)
\\
&&\quad\qquad{}+ \mathrm{E}\bigl\vert\mathbf{x}_n^*\beta_1(z)
\mathbf{A}_1^{-1}(z)\mathbf{r}_1
\mathbf{r}_1^*\mathbf{A}_1^{-1}(z)
\mathbf{x}_n\bigr\vert^{2l} \mathrm{I}\bigl(\bigl|
\beta_1(z)\bigr|>C\bigr)
\\
&&\qquad \leq \frac{C}{N^{2l}}\mathrm{E}\bigl\vert\mathbf{X}_1^*\mathbf
{A}_1^{-1}(z)\mathbf{x}_n\mathbf{x}_n^*
\mathbf{A}_1^{-1}(z)\mathbf{X}_1\bigr\vert
^{2l}+o\bigl(N^{-t}\bigr)
\\
&&\qquad \leq \frac{C}{N^{2l}}\mathrm{E} \bigl\vert\mathbf{X}_1^*
\mathbf{A}_1^{-1}(z)\mathbf{x}_n
\mathbf{x}_n^*\mathbf{A}_1^{-1}(z)\mathbf
{X}_1-\mathbf{x}_n^*\mathbf{A}_1^{-2}(z)
\mathbf{x}_n\bigr\vert^{2l}\\
&&\quad\qquad{} +\frac
{C}{N^{2l}}\mathrm{E}
\bigl\vert\mathbf{x}_n^*\mathbf{A}_1^{-2}(z)
\mathbf{x}_n\bigr\vert^{2l}
\\
&&\qquad \leq\frac{C}{N^{2l}} \bigl[ \mathrm{E} \bigl(\operatorname{tr}
\mathbf{A}_1^{-1}(z)\mathbf{x}_n\mathbf{x}_n^*
\mathbf{A}_1^{-1}(z)\mathbf{A}_1^{-1}(
\bar{z})\mathbf{x}_n\mathbf{x}_n^*\mathbf{A}_1^{-1}(
\bar{z}) \bigr)^{l}
\\
&&\quad\qquad\qquad{}+N^{l-2}\mathrm{E}\operatorname{tr} \bigl( \mathbf{A}_1^{-1}(z)
\mathbf{x}_n\mathbf{x}_n^*\mathbf{A}_1^{-1}(z)
\mathbf{A}_1^{-1}(\bar{z})\mathbf{x}_n\mathbf
{x}_n^*\mathbf{A}_1^{-1}(\bar{z})
\bigr)^l \bigr]
\\
&&\quad\qquad{} +\frac{C}{N^{2l}}\mathrm{E}\bigl\vert
\mathbf{x}_n^*\mathbf{A}_1^{-2}(z)
\mathbf{x}_n\bigr\vert^{2l}
\\
&&\qquad \leq\frac{C}{N^{l+2}v^{2l}}\mathrm{E}\bigl\vert m_{H^{\mathbf
{B}_1}}(z)\bigr\vert^{2l}.
\end{eqnarray*}
The last step follows the fact that
\[
\mathbf{x}_n^*\mathbf{A}_1^{-1}(z)
\mathbf{A}_1^{-1}(\bar{z})\mathbf{x}_n=v^{-1}
\Im\bigl(\mathbf{x}_n^*\mathbf{A}_1^{-1}(z)
\mathbf{x}_n\bigr) =v^{-1}\Im\bigl(m_{H^{\mathbf{B}_1}}(z)
\bigr).
\]

For $v^2v_y\geq C_0N^{-1}$, which implies that $v\geq C_0N^{-1/2}$.
Choose $C_0$ and $N$ large enough, such that $\frac
{C}{N^{l+2}v^{2l}}\leq\frac{1}{2}$. Further, by $C_r$-inequality, we
obtain
\begin{eqnarray*}
\mathrm{E}\bigl\vert m_{H^{\mathbf{B}_1}}(z)\bigr\vert^{2l}&\leq& C
\mathrm{E}\bigl\vert m_{H^{\mathbf{B}_1}}(z)-m_n^H(z)\bigr
\vert^{2l} +C\mathrm{E}\bigl\vert m_n^H(z)
\bigr\vert^{2l}
\\
&\leq&\tfrac{1}{2}\mathrm{E}\bigl\vert m_{H^{\mathbf{B}_1}}(z)\bigr
\vert^{2l}+C\mathrm{E}\bigl\vert m_n^H(z)\bigr
\vert^{2l}.
\end{eqnarray*}
That is, $\mathrm{E}\vert m_{H^{\mathbf{B}_1}}(z)\vert^{2l} \leq
C\mathrm{E}\vert m_n^H(z) \vert^{2l}$, for some constant $C$. This
finishes the proof.
\end{pf}

%
%
\begin{lem}\label{Ap5}
Under the conditions of Theorem \ref{theorem2}, for $v^2v_y\geq
C_0N^{-1}$, we have
\[
\mathrm{E}\bigl\vert m_n^H(z)-\mathrm{E}m_n^H(z)\bigr\vert^{2l} \leq
\frac{C}{N^lv^{2l}} \biggl(\frac1{v_y}+\frac{\Delta^H}{v}
\biggr)^{2l}.
\]
\end{lem}

\begin{pf}
Write $\mathrm{E}_j(\cdot)$ as the conditional expectation given
$\{\mathbf{r}_{1},\ldots,\mathbf{r}_{j} \}$. It can then be shown that
$ m_n^H(z)-\mathrm{E}m_n^H(z)=\sum_{j=1}^N \gamma_j, $ where
\begin{eqnarray*}
\gamma_j&=:& \mathrm{E}_j \bigl(\mathbf{x}_n^*
\mathbf{A}^{-1}(z)\mathbf{x}_n \bigr)-\mbox
{E}_{j-1} \bigl(\mathbf{x}_n^*\mathbf{A}^{-1}(z)
\mathbf{x}_n \bigr)
\\
&=& (\mathrm{E}_j-\mathrm{E}_{j-1} ) \bigl\{
\mathbf{x}_n^* \bigl(\mathbf{A}^{-1}(z)-
\mathbf{A}_j^{-1}(z) \bigr)\mathbf{x}_n \bigr
\}
\\
&=&- (\mathrm{E}_j-\mathrm{E}_{j-1} ) \bigl\{
\beta_j(z) \mathbf{x}_n^*\mathbf{A}_j^{-1}(z)
\mathbf{r}_j \mathbf{r}_j^*\mathbf{A}_j^{-1}(z)
\mathbf{x}_n \bigr\}.
\end{eqnarray*}
Therefore, by Lemmas \ref{B3}(b), we have
\[
\mathrm{E}\bigl\vert m_n^H(z)-\mathrm{E}m_n^H(z)
\bigr\vert^{2l} \leq C\mathrm{E}\Biggl(\sum_{j=1}^N
\mathrm{E}_{j-1}|\gamma_j|^2
\Biggr)^l+C\sum_{j=1}^N\mathrm{E}|
\gamma_j|^{2l}.
\]
Using Lemmas \ref{Ap3} and \ref{B4}, we have
\begin{eqnarray*}
\mathrm{E}_{j-1}&=&\mathrm{E}_{j-1} \bigl\vert(
\mathrm{E}_j-\mathrm{E}_{j-1}) \beta_j(z)
\mathbf{x}_n^*\mathbf{A}_j^{-1}(z)\mathbf
{r}_j\mathbf{r}_j^*\mathbf{A}_j^{-1}(z)
\mathbf{x}_n\bigr\vert^2
\\
&\leq&\frac{C}{N^2}\mathrm{E}_{j-1}\bigl\vert
\mathbf{X}_j^*\mathbf{A}_j^{-1}(z)
\mathbf{x}_n\mathbf{x}_n^*\mathbf{A}_j^{-1}(z)
\mathbf{X}_j -\mathbf{x}_n^*\mathbf{A}_j^{-2}(z)
\mathbf{x}_n\bigr\vert^2
\\
&&{} +\frac{C}{N^2}\mathrm{E}_{j-1}\bigl\vert\mathbf{x}_n^*
\mathbf{A}_j^{-2}(z)\mathbf{x}_n\bigr\vert^2
\\
&\leq&\frac{C}{N^2}\mathrm{E}_{j-1}\operatorname{tr} \bigl(\mathbf
{A}_j^{-1}(z)\mathbf{x}_n
\mathbf{x}_n^*\mathbf{A}_j^{-1}(z)\mathbf
{A}_j^{-1}(\bar{z})\mathbf{x}_n
\mathbf{x}_n^*\mathbf{A}_j^{-1}(\bar{z})
\bigr)
\\
&&{} +\frac{C}{N^2}\mathrm{E}_{j-1}\bigl\vert
\mathbf{x}_n^*\mathbf{A}_j^{-2}(z)
\mathbf{x}_n\bigr\vert^2.
\end{eqnarray*}
By the fact that $\mathbf{x}_n^*\mathbf{A}_j^{-1}(z)\mathbf
{A}_j^{-1}(\bar{z})\mathbf{x}_n=v^{-1}\Im(\mathbf{x}_n^*\mathbf
{A}_j^{-1}(z)\mathbf{x}_n)$ and $\|\mathbf{A}_j^{-1}(z)\|\leq v^{-1}$,
we have
\[
\mathrm{E}_{j-1}|\gamma_j|^2\leq
\frac{C}{N^2v^2}\mathrm{E}_{j-1}\bigl\vert m_{H^{\mathbf
{B}_j}}(z)\bigr
\vert^2.
\]
On the other side,
\begin{eqnarray*}
\mathrm{E}|\gamma_j|^{2l}&=&\frac{1}{N^{2l}}\mathrm{E}
\bigl\vert\mathbf{X}_1^*\mathbf{A}_1^{-1}(z)
\mathbf{x}_n\mathbf{x}_n^*\mathbf{A}_1^{-1}(z)
\mathbf{X}_1\bigr\vert^{2l}
\\
&\leq&\frac{C}{N^{2l}}\mathrm{E} \bigl\vert\mathbf{X}_1^*
\mathbf{A}_1^{-1}(z)\mathbf{x}_n
\mathbf{x}_n^*\mathbf{A}_1^{-1}(z)\mathbf
{X}_1-\mathbf{x}_n^*\mathbf{A}_1^{-2}(z)
\mathbf{x}_n\bigr\vert^{2l}
\\
&&{}+\frac
{C}{N^{2l}}\mathrm{E}
\bigl\vert\mathbf{x}_n^*\mathbf{A}_1^{-2}(z)
\mathbf{x}_n\bigr\vert^{2l}
\\
&\leq&\frac{C}{N^{2l}}\times\frac{N^{l-2}}{v^{2l}}\mathrm{E}
\bigl(\Im\bigl(
\mathbf{x}_n^*\mathbf{A}_1^{-1}(z)
\mathbf{x}_n\bigr) \bigr)^{2l} +\frac
{C}{N^{2l}}\mathrm{E}
\bigl\vert\mathbf{x}_n^*\mathbf{A}_1^{-1}(z)
\mathbf{x}_n\bigr\vert^{2l}
\\
&\leq&\frac{C}{N^{l+2}v^{2l}}\mathrm{E}\bigl\vert m_{H^{\mathbf
{B}_1}}(z)\bigr\vert
^{2l}.
\end{eqnarray*}
Thus, we obtain
\begin{eqnarray*}
\mathrm{E}\bigl\vert m_n^H(z)-\mathrm{E}m_n^H(z)
\bigr\vert^{2l}&\leq& \frac
{C}{N^lv^{2l}}\mathrm{E}\bigl\vert
m_{H^{\mathbf{B}_1}}(z)\bigr\vert^{2l} +\frac
{C}{N^{l+1}v^{2l}}\mathrm{E}
\bigl\vert m_{H^{\mathbf{B}_1}}(z)\bigr\vert^{2l}
\\
&\leq&\frac{C}{N^lv^{2l}}\mathrm{E}\bigl\vert m_n^H(z)
\bigr\vert^{2l}\qquad\mbox{(by Lemma \ref{Ap4})}.
\end{eqnarray*}
Further
\[
\mathrm{E}\bigl\vert m_n^H(z)\bigr\vert
^{2l}\leq\mathrm{E}\bigl\vert m_n^H(z)-\mathrm{E}m_n^H(z)\bigr\vert^{2l}+\bigl\vert
\mathrm{E}m_n^H(z)-m_y(z)\bigr\vert
^{2l}+\bigl\vert m_y(z)\bigr\vert^{2l}.
\]
For $v\geq C_0N^{-1/2}$, choose $C_0$ large enough, such that $\frac
{C}{N^lv^{2l}}\leq\frac{1}{2}$. And using integration by parts, it is
easy to find that
%
%
\begin{equation}
\label{eqAp2}
\bigl\vert\mathrm{E}m_n^H(z)-m_y(z)\bigr \vert\leq\frac{C\Delta^H}{v},
\end{equation}
where $\Delta^H=\|EH^{\mathbf{S}_n}-F_{y_n}\|$.

Besides, from Lemma \ref{B7}, we know that $\vert m_y(z)\vert\leq
\frac{C}{v_y}$.\eject

Therefore, we obtain
\[
\mathrm{E}\bigl\vert m_n^H(z)-\mathrm{E}m_n^H(z)\bigr\vert^{2l} \leq
\frac{C}{N^lv^{2l}} \biggl(\frac1{v_y}+\frac{\Delta^H}{v}
\biggr)^{2l}.
\]
The proof is then complete.
\end{pf}

%
%
\begin{lem}\label{Ap6}
\[
\mathrm{E} \bigl|\alpha_1(z) \bigr|^2 \leq\frac{C}{N^2v}
\biggl(\frac1{v_y}+\frac{\Delta^H}{v} \biggr).
\]
\end{lem}

\begin{pf}
Lemma \ref{B4} implies that
\begin{eqnarray*}
\mathrm{E} \bigl|\alpha_1(z) \bigr|^2&=& \frac{1}{N^2}
\mathrm{E}\bigl\vert\mathbf{X}_1^*\mathbf{A}_1^{-1}(z)
\mathbf{x}_n \mathbf{x}_n^* \mathbf{X}_1 -
\mathbf{x}_n^*\mathbf{A}_1^{-1}(z)
\mathbf{x}_n\bigr\vert^2
\\
&\leq& \frac{C}{N^2}\mathrm{E} \bigl(\mathbf{x}_n^*
\mathbf{A}_1^{-1}(\bar{z}) \mathbf{A}_1^{-1}(z)
\mathbf{x}_n \bigr)
\\
&\leq& \frac{C}{N^2v}\bigl\vert\mathrm{E}m_{H^{\mathbf
{B}_1}}(z)\bigr\vert.
\end{eqnarray*}
Using Lemmas \ref{Ap4} and \ref{Ap5} and integration by parts, we have
\begin{eqnarray*}
\mathrm{E}\bigl\vert m_{H^{\mathbf{B}_1}}(z)\bigr\vert^2&\leq&
\mathrm{E}\bigl\vert m_n^H(z)\bigr\vert^2
\\
&\leq&\mathrm{E}\bigl\vert m_n^H(z)-
\mathrm{E}m_n^H(z)\bigr\vert^2+\bigl\vert
\mathrm{E}m_n^H(z)-m_y(z)\bigr\vert
^2+\bigl\vert m_y(z)\bigr\vert^2
\\
&\leq& C \biggl(\frac{1}{v_y}+\frac{\Delta^H}{v} \biggr)^2.
\end{eqnarray*}
This finishes the proof.
\end{pf}

%
%
\begin{lem}\label{Ap7}
Under the conditions in Theorem \ref{theorem1}, for any fixed $t>0$,
\[
\int_B^{\infty} \bigl\vert\mathrm{E}H^{\mathbf{S}_n}(x)-F_{y_n}(x)\bigr
\vert\,dx=o\bigl(N^{-t}
\bigr).
\]
\end{lem}

\begin{pf}
For any fixed $t>0$, by Lemma \ref{913}, it follows that
\[
\mathrm{P} \bigl(\lambda_{\max}(\mathbf{S}_n) \geq B+x \bigr)
\leq C N^{-t-1}(B+x-\varepsilon)^{-2}
\]
and
\begin{eqnarray*}
\int_B^{\infty} \bigl\vert\mathrm{E}H^{\mathbf{S}_n}(x)-F_{y_n}(x)
\bigr\vert\,dx &\leq&\int_B^{\infty} \bigl(1-
\mathrm{E}H^{\mathbf{S}_n}(x) \bigr)\,dx
\\
&=&\int_B^{\infty} \Biggl(1-\sum
_{i=1}^N |y_i|^2\mathrm{P}(
\lambda_i\leq x) \Biggr)\,dx
\\
&\leq& \int_B^{\infty} \Biggl( \sum
_{i=1}^N |y_i|^2-\sum
_{i=1}^N |y_i|^2
\mathrm{P}(\lambda_{\max}\leq x) \Biggr)\,dx
\\
&\leq& \int_B^{\infty} N^{-t-1}(B+x-
\varepsilon)^{-2}\,dx=o\bigl(N^{-t}\bigr).
\end{eqnarray*}
The proof is complete.
\end{pf}

\section{}\label{sec8}
In what follows, we will present some existing results which are of
substantial importance in proving the main theorems.

%
%
\begin{lem}[(Theorem 2.2 in \cite{bai1993a})]\label{LL}
Let $F$ be a distribution function and let $G$ be a function of bounded
variation satisfying $\int|F(x)-G(x)|\,dx < \infty$. Denote their
Stieltjes transforms by $f(z)$ and $g(z)$, respectively. Then
\begin{eqnarray*}
\|F-G\|&=& \sup_x \bigl|F(x)-G(x)\bigr|
\\
&\leq& \frac{1}{\pi(1-\kappa)(2\gamma-1)}
\\
&&{}\times \biggl(\int_{-A}^{A}\bigl|f(z)-g(z)\bigr|\,du +\frac{2\pi}{v}\int_{|x|>B} \bigl|F(x)-G(x)\bigr|\,dx
\\
&&\hspace*{76pt}{}+\frac{1}{v}\sup_x \int_{|y|\leq2v\tau}
\bigl|G(x+y)-G(x)\bigr|\,dy \biggr),
\end{eqnarray*}
where $z=u+iv$ is a complex variable, $\gamma$, $\kappa$, $\tau$, $A$
and $B$ are positive constants such that $A>B$,
\[
\kappa=\frac{4B}{\pi(A-B)(2\gamma-1)}<1
\]
and
\[
\gamma=\frac{1}{\pi}
\int_{|u|<\tau}\frac{1}{u^2+1}\,du>\frac{1}{2}.
\]
\end{lem}

%
%
\begin{lem}[(Lemma 8.15 in \cite{baijack2010})]\label{B1}
For any $v>0$, we have
\begin{eqnarray*}
\sup_x \int_{|u|<v} \bigl|F_{y_n}(x+u)-F_{y_n}(x)\bigr|\,du &<&\frac{11\sqrt{2(1+y_n)}}{3\pi y_n}v^2/v_y,
\end{eqnarray*}
where $F_{y_n}$ is the c.d.f. of the Mar\v{c}enko--Pastur distribution
with index $y_n\leq1$, and $v_y=1-\sqrt{y_n}+\sqrt{v}$.
\end{lem}

%
%
\begin{lem}[((1.15) in \cite{baijack2004})]\label{B2}
Let $\mathbf{A}=(a_{ij})_{n\times n}$ and $\mathbf
{B}=(b_{ij})_{n\times n}$ be two nonrandom matrices. Let
$\mathbf{X}=(X_1,\ldots,X_n)^*$ be a random vector of independent
complex entries. Assume that $\mathrm{E}X_i=0$ and
$\mathrm{E}|X_i|^2=1$. Then we have
\begin{eqnarray*}
&& \mathrm{E}\bigl(\mathbf{X}^*\mathbf{A}\mathbf{X}-\operatorname{tr}
\mathbf{A}\bigr)
\bigl(\mathbf{X}^*\mathbf{B}\mathbf{X}-\operatorname{tr} \mathbf{B}\bigr)
\\
&&\qquad = \sum_{i=1}^n \bigl(\mathrm{E}|X_i|^4-|
\mathrm{E}X_i^2|^2-2\bigr)a_{ii}b_{ii}
+\bigl|\mathrm{E}X_i^2\bigr|^2\operatorname{tr}
\mathbf{A} \mathbf{B}^T+\operatorname{tr} \mathbf{AB}.
\end{eqnarray*}
\end{lem}

%
%
\begin{lem}[(Burkholder inequalities (Lemmas 2.1 and 2.2 in \cite{bai1998}))]\label{B3} Let $\{X_k\}$ be a complex martingale
difference sequence with respect to the increasing $\sigma$-field
$\{\mathcal{F}_k\}$, and let $\mathrm{E}_k$ denote the conditional
expectation with respect to $\mathcal{F}_k$. Then we have:
\begin{longlist}[(a)]
\item[(a)] for $p>1$,
\[
\mathrm{E} \Biggl|\sum_{k=1}^n
X_k \Biggr|^p \leq K_p \mathrm{E} \Biggl(\sum
_{k=1}^n |X_k|^2
\Biggr)^{p/2},
\]

\item[(b)] for $p \geq2$,
\[
\mathrm{E} \Biggl|\sum_{k=1}^n
X_k \Biggr|^p \leq K_p \Biggl(\mathrm{E}
\Biggl(\sum_{k=1}^n \mathrm{E}_{k-1}|X_k|^2 \Biggr)^{p/2}+\mbox{
\em E}\sum_{k=1}^n |X_k|^p
\Biggr),
\]
where $K_p$ is a constant which depends on $p$ only.
\end{longlist}
\end{lem}

%
%
\begin{lem}[(Lemma 2.7 in \cite{bai1998})]\label{B4}
Let $\mathbf{A}=(a_{ij})$ be an $n\times n$ nonrandom matrix and
$\mathbf{X}=(X_1,\ldots, X_n)^*$ be random vector of independent
complex entries. Assume that $\mathrm{E}X_i=0$, $\mathrm{E}|X_i|^2=1$
and $\mathrm{E}|X_i|^l \leq V_l$. Then for any $p\geq2$,
\[
\mathrm{E}\bigl|\mathbf{X}^*\mathbf{A}\mathbf{X}-\operatorname{tr}
\mathbf
{A}\bigr|^p \leq K_p \bigl( \bigl(V_4\operatorname{tr}\bigl(\mathbf{AA}^*\bigr)
\bigr)^{p/2}+V_{2p}\operatorname{tr}
\bigl(\mathbf{A}\mathbf{A}^*\bigr)^{p/2} \bigr),
\]
where $K_p$ is a constant depending on $p$ only.
\end{lem}

%
%
\begin{lem}[(Lemma 2.6 in \cite{jackbai1995})]\label{B6}
Let $z\in\mathbb{C}^+$ with $v=\Im(z)$, $\mathbf{A}$ and $\mathbf{B}$
$n\times n$ with $\mathbf{B}$ Hermitian, $\tau\in\mathbb{R}$,
and $\mathbf{q} \in\mathbb{C}^n$. Then
\[
\bigl| \operatorname{tr}\bigl((\mathbf{B}-z\mathbf{I}_n)^{-1}-
\bigl(\mathbf{B}+\tau\mathbf{q}\mathbf{q}^*-z\mathbf{I}_n
\bigr)^{-1}\bigr)\mathbf{A} \bigr| \leq\frac{\|\mathbf{A}\|}{v},
\]
where $\|A\|$ denotes spectral norm on matrices.
\end{lem}

%
%
\begin{lem}[((8.4.9) in \cite{baijack2010})]\label{B7}
For the Stieltjes transform of the Mar\v{c}enko--Pastur distribution,
we have
\[
\bigl|m_y(z)\bigr| \leq\frac{\sqrt{2}}{\sqrt{y}v_y},
\]
where $v_y=\sqrt{a}+\sqrt{v}=1-\sqrt{y_n}+\sqrt{v}$.
\end{lem}

%
%
\begin{lem}[(Lemma 8.20 in \cite{baijack2010})]\label{B8}
If $|z|<A$, $v^2v_y\geq C_0N^{-1}$ and $l \geq1$, then
\[
\mathrm{E}\bigl|m_n(z)-\mathrm{E}m_n(z)\bigr|^{2l}
\leq\frac{C}{N^{2l}v^{4l}y_n^{2l}} \biggl(\Delta+\frac{v}{v_y}
\biggr)^l,
\]
where A is a positive constant, $v_y=1-\sqrt{y_n}+\sqrt{v}$ and
$\Delta:=\|\mathrm{E}F^{\mathbf{S}_n}-F_{y_n} \|$.
\end{lem}

%
%
\begin{lem}[(Lemma 9.1 in \cite{baijack2010})]\label{B10}
Suppose that $X_i$, $i=1,\ldots, n$, are independent, with
$\mathrm{E}X_i=0$, $\mathrm{E}|X_i|^2=1$, $\sup\mathrm{E}|X_i|^4=\nu<
\infty$ and $|X_i|\leq\eta\sqrt{n}$ with $\eta>0$. Assume that
$\mathbf{A}$ is a complex matrix. Then for any given $p$ such that
$2\leq p \leq b \log(n\nu^{-1}\eta^4 )$ and $b>1$, we have
\[
\mathrm{E}\bigl|\bolds{\alpha}^* \mathbf{A} \bolds{\alpha} -
\operatorname{tr} (\mathbf{A})\bigr|^p \leq\nu n^p\bigl(n\eta^4
\bigr)^{-1} \bigl(40b^2\|\mathbf{A}\|\eta^2
\bigr)^p,
\]
where $\bolds{\alpha}=(X_1,\ldots, X_n)^T$.
\end{lem}

%
%
\begin{lem}[(Theorem 8.10 in \cite{baijack2010})]\label{B12}
Let $\mathbf{S}_n=\mathbf{X}\mathbf{X}^*/N$, where $\mathbf {X}=
(X_{ij}(n))_{n\times N}$. Assume that the following conditions hold:
\begin{longlist}[(3)]
\item[(1)] For each $n$, $X_{ij} (n)$ are independent,

\item[(2)] $\mathrm{E}X_{ij}(n)=0$, $\mathrm{E}|X_{ij}(n)|^2=1$, for
    all $i,j$,

\item[(3)] $\sup_n\sup_{i,j} \mathrm{E}|X_{ij}(n)|^6<\infty$.
\end{longlist}
Then we have
\[
\Delta=: \bigl\|\mathrm{E}F^{\mathbf{S}_n}-F_{y_n} \bigr\|= \cases{ O
\bigl(N^{-1/2}a^{-1}\bigr), &\quad if $a>N^{-1/3}$,
\vspace*{2pt}
\cr
O\bigl(N^{-1/6}\bigr), &\quad otherwise,}
\]
where $y_n=n/N\leq1$ and $a$ is defined in the Mar\v{c}enko--Pastur
distribution.
\end{lem}

%
%
\begin{lem}[(Theorem 5.11 in \cite{baijack2010})]\label{B14}
Assume that the entries of $\{X_{ij}\}$ is a double array of i.i.d.
complex random variables with mean zero, variance $\sigma^2$ and finite
4th moment. Let $\mathbf{X}=(X_{ij})_{n\times N}$ be the $n\times N$
matrix of the upper-left corner of the double array. If $n/N\rightarrow
y \in(0,1)$, then, with probability one, we have
\[
\lim_{n\rightarrow\infty} \lambda_{\min}(\mathbf{S}_n)=
\sigma^2(1-\sqrt{y})^2
\]
and
\[
\lim_{n\rightarrow\infty} \lambda_{\max}(\mathbf{S}_n)=
\sigma^2(1+\sqrt{y})^2.
\]
\end{lem}

%
%
\begin{lem}[(Theorem 5.9 in \cite{baijack2010})]\label{913}
Suppose that the entries of the matrix $\mathbf{X}=(X_{ij})_{n\times
N}$ are independent (not necessarily identically distributed) and
satisfy:
\begin{longlist}[(3)]
\item[(1)] $\mathrm{E}X_{ij}=0$,

\item[(2)] $|X_{ij}|\leq\sqrt{N}\delta_N$,\eject

\item[(3)] $\max_{ij}|\mathrm{E}|X_{ij}|^2-\sigma^2|\rightarrow0$ as
    $N\rightarrow\infty$ and

\item[(4)] $\mathrm{E}|X_{ij}|^l\leq b(\sqrt{N}\delta_N)^{l-3}$ for
    all $l\geq3$,
    where $\delta_N\rightarrow0$ and $b>0$. Let $\mathbf{S}_n=\mathbf
{X}\mathbf{X}^*/N$. Then, for any
$x>\epsilon>0$, $n/N\rightarrow y$, and fixed integer $\ell\geq2$,
we have
\[
\mathrm{P} \bigl(\lambda_{\max}(\mathbf{S}_n)\geq
\sigma^2(1+\sqrt{y})^2+x \bigr)\leq CN^{-\ell}
\bigl(\sigma^2(1+\sqrt{y})^2+x-\epsilon
\bigr)^{-\ell}
\]
for some constant $C>0$.
\end{longlist}
\end{lem}

\section{}\label{sec9}
Note that the data matrix $\mathbf{X}=(X_{ij})_{n\times N}$ consists of
i.i.d. complex random variables with mean 0 and variance 1. In what
follows, we will further assume that every $|X_{ij}|$ is bounded by
$\eta_NN^{1/4}$ for some carefully selected $\eta_N$. The proofs
presented in the following three steps jointly justify such a
convenient assumption.

\subsection{\texorpdfstring{Truncation for Theorem \protect\ref{theorem1}}{Truncation for Theorem 1.1}}\label{sec9.1}
Choose $\eta_N \downarrow0$ and $\eta_NN^{1/4}\uparrow\infty$ as
$N\to\infty$ such that
\[
\lim_{N \rightarrow\infty}\eta_N^{-10}\mathrm{E}|X_{11}|^{10}\mathrm{I} \bigl(|X_{11}|>
\eta_NN^{1/4} \bigr)=0.
\]
Let $\widehat{\underline{\mathbf{X}}}_n$ denote the truncated data
matrix whose entry on the $i$th row and $j$th column is
$X_{ij}\mathrm{I} (|X_{ij}|\leq\eta_NN^{1/4} )$, $i=1,\ldots,n$,
$j=1,\ldots, N$. Define
$\widehat{\mathbf{S}}_n=\widehat{\underline{\mathbf{X}}}_n\widehat
{\underline{\mathbf{X}}}_n^*/N$. Then
\begin{eqnarray*}
\mathrm{P} (\mathbf{S}_n\neq\widehat{\mathbf{S}}_n ) &\leq&
nN \mathrm{P} \bigl(|X_{ij}|>\eta_NN^{1/4} \bigr)
\\
&\leq& nN^{-3/2}\eta_N^{-10}\mathrm{E}|X_{11}|^{10}\mathrm{I} \bigl(|X_{11}|>
\eta_NN^{1/4} \bigr)
\\
&=&o\bigl(N^{-1/2}\bigr).
\end{eqnarray*}

\subsection{\texorpdfstring{Truncation for Theorems \protect\ref{theorem2} and
\protect\ref{theorem3}}{Truncation for Theorems 1.6 and 1.8}}\label{sec9.2}
Choose $\eta_N
\downarrow0$ and $\eta_NN^{1/4}\uparrow\infty$ as $N\to\infty$
such that
%
%
\begin{eqnarray}\label{eta2}
\lim_{N \rightarrow\infty}\eta_N^{-8}\mathrm{E}|X_{11}|^{8}\mathrm{I} \bigl(|X_{11}|>
\eta_NN^{1/4} \bigr)&=&0.
\end{eqnarray}
Let $\widehat{\underline{\mathbf{X}}}_n$ denote the truncated data
matrix whose entry on the $i$th row and $j$th column is
$X_{ij}\mathrm{I} (|X_{ij}|\leq\eta_NN^{1/4} )$, $i=1,\ldots,n$,
$j=1,\ldots, N$. Define
$\widehat{\mathbf{S}}_n=\widehat{\underline{\mathbf{X}}}_n\widehat
{\underline{\mathbf{X}}}_n^*/N$. Then
\begin{eqnarray*}
\mathrm{P} (\mathbf{S}_n\neq\widehat{\mathbf{S}}_n\mbox{, i.o.}) &=&\lim_{k\rightarrow\infty} \mathrm{P} \Biggl(\bigcup
_{N=k}^{\infty} \bigcup_{i=1}^n
\bigcup_{j=1}^N |X_{ij}|>
\eta_NN^{1/4} \Biggr)
\\
&=&\lim_{k\rightarrow\infty} \mathrm{P} \Biggl(\bigcup
_{t=k}^{\infty} \bigcup_{N\in[2^t,2^{t+1})}
\bigcup_{i=1}^n \bigcup
_{j=1}^N|X_{ij}|>\eta_NN^{1/4}
\Biggr)
\\
&\leq& \lim_{k\rightarrow\infty} \sum_{t=k}^{\infty}
\mathrm{P} \Biggl(\bigcup_{i=1}^{(y_n+1) 2^{t+1}} \bigcup
_{j=1}^{2^{t+1}} |X_{ij}|>
\eta_{2^t}2^{t/4} \Biggr)
\\
&\leq& C\lim_{k\rightarrow\infty} \sum_{t=k}^{\infty}
\bigl(2^{t+1}\bigr)^2 \mathrm{P} \bigl(|X_{11}|>
\eta_{2^t}2^{t/4} \bigr)
\\
&\leq& C\lim_{k\rightarrow\infty} \sum_{t=k}^{\infty}
\sum_{l=t}^{\infty} 4^t \mathrm{P}
\bigl(\eta_{2^l}2^{l/4}<|X_{11}|\leq\eta
_{2^{l+1}}2^{(l+1)/4} \bigr)
\\
&=&C\lim_{k\rightarrow\infty} \sum_{l=k}^{\infty}
\sum_{t=k}^{l} 4^t \mathrm{P}
\bigl(\eta_{2^l}2^{l/4}<|X_{11}|\leq\eta
_{2^{l+1}}2^{(l+1)/4} \bigr)
\\
&\leq& \lim_{k\rightarrow\infty} \sum_{l=k}^{\infty}C
\eta_{2^l}^{-8} \mathrm{E}|X_{11}|^8
\mathrm{I} \bigl(\eta_{2^l}2^{l/4}<|X_{11}|\leq
\eta_{2^{l+1}}2^{(l+1)/4} \bigr)
\\
&=&0.
\end{eqnarray*}
The last equality is due to (\ref{eta2}).

\subsection{Centralization}\label{sec9.3}
The centralization procedures for three theorems are identical, only
8th moment is required and thus
we treat them uniformly. Let $\widetilde{\underline{\mathbf{X}}}_n$
denote the centralized version of
$\widehat{\underline{\mathbf{X}}}_n$. More explicitly, on the $i$th
row and $j$th column of $\widetilde{\underline{\mathbf{X}}}_n$, the
entry is
\[
X_{ij}\mathrm{I} \bigl(|X_{ij}|\leq\eta_N
N^{1/4} \bigr)-\mathrm{E} \bigl(X_{ij}\mathrm{I}
\bigl(|X_{ij}|\leq\eta_N N^{1/4} \bigr) \bigr).
\]

Notice that according to Theorem 3.1 of \cite{ybk1988},
$\|(\mathbf{S}_n-z\mathbf{I}_n)^{-1} \|$ is bounded by~$1/v$, where $\|
\cdot\|$ denotes the spectral norm for a matrix. Define
$\widetilde{\mathbf{S}}_n=\widetilde{\underline{\mathbf{X}}}_n
\widetilde{\underline{\mathbf{X}}}{}_n^*/N$. Suppose that $v\geq
C_0N^{-1/2}$, we obtain
\begin{eqnarray*}
&&\bigl\vert m_{H^{\widehat{\mathbf{S}}_n}}(z)-m_{H^{\widetilde
{\mathbf{S}}_n}}(z)\bigr\vert
\\
&&\qquad =\bigl\vert\mathbf{x}_n^* (\widehat{\mathbf
{S}}_n-z\mathbf{I}_n )^{-1}
\mathbf{x}_n- \mathbf{x}_n^* (\widetilde{
\mathbf{S}}_n-z\mathbf{I}_n )^{-1}
\mathbf{x}_n\bigr\vert
\\
&&\qquad \leq \bigl\llVert(\widehat{\mathbf{S}}_n-z\mathbf{I}_n)^{-1}
\bigr\rrVert\llVert\widehat{\mathbf{S}}_n-\widetilde{
\mathbf{S}}_n\rrVert\bigl\llVert(\widetilde{\mathbf{S}}_n-z
\mathbf{I}_n)^{-1}\bigr\rrVert
\\
&&\qquad \leq \frac{1}{v^2}\llVert\widehat{\mathbf{S}}_n-\widetilde{
\mathbf{S}}_n\rrVert
\\
&&\qquad \leq \frac{1}{Nv^2} \bigl(\llVert\widehat{\mathbf{\underline
{X}}}_n\rrVert\bigl\llVert\widehat{\underline{\mathbf{X}}}{}_n^*-
\widetilde{\underline{\mathbf{X}}}{}_n^*\bigr\rrVert+\llVert
\widehat{\underline{\mathbf{X}}}_n-\widetilde{\mathbf{\underline
{X}}}_n\rrVert\bigl\llVert\widetilde{\underline{\mathbf{X}}}{}_n^*
\bigr\rrVert\bigr)\qquad\mbox{by Lemma \ref{913}}
\\
&&\qquad \leq \frac{C}{\sqrt{N}v^2} \llVert\widehat{\underline{\mathbf{X}}}_n-
\widetilde{\underline{\mathbf{X}}}_n\rrVert\qquad\mbox{a.s.}
\\
&&\qquad =\frac{C}{\sqrt{N}v^2}\bigl\vert\mathrm{E} \bigl\{X_{11} \mathrm{I}
\bigl(|X_{11}|\leq\eta_NN^{1/4} \bigr) \bigr\}
\bigr\vert\llVert\mathbf{1}_{n\times1}\rrVert\bigl\llVert
\mathbf{1}_{N\times1}^{\prime}\bigr\rrVert
\\
&&\qquad \leq C\sqrt{N}v^{-2}\eta_N^{-7}
N^{-7/4} \mathrm{E} \bigl( |X_{11}|^8 \mathrm{I}
\bigl(|X_{11}|> \eta_NN^{1/4} \bigr) \bigr)
\\
&&\qquad = o\bigl(N^{-1/4}\bigr).
\end{eqnarray*}

To establish both the weak and the strong convergence rates of the VESD
to the Mar\v{c}enko--Pastur distribution, this $o(N^{-1/4})$ suffices.
Moreover, for the convergence rate presented in Theorem \ref{theorem1},
we shall prove the following. Let $m_y(z)$ denotes the Stieltjes
transform of the Mar\v{c}enko--Pastur distribution, thus $m_y(z)$ is
bounded by $\frac{\sqrt{2}}{\sqrt{y_n}(1-\sqrt{y_n}+\sqrt{v})}\leq
\frac{C}{\sqrt{v}}$ in Lemma \ref{B7}. Then
\[
\bigl|\mathrm{E}m_n^H(z)\bigr|\leq\bigl|\mathrm{E}m_n^H(z)-m_y(z)\bigr|+\bigl|m_y(z)\bigr|
\leq C\bigl|m_y(z)\bigr|\leq\frac{C}{\sqrt{v}}.
\]
Besides
$\mathbf{x}_n^*(\widehat{\mathbf{S}}_n-z\mathbf{I}_n)^{-1}\mathbf{x}_n$
can be considered as a Stieltjes transform of some VESD function. So,
we have
\[
\mathrm{E}\bigl\|\mathbf{x}_n^*(\widehat{\mathbf{S}}_n-z
\mathbf{I}_n)^{-1}\bigr\|^2 =v^{-1}
\mathrm{E}\mathbf{x}_n^*(\widehat{\mathbf{S}}_n-z\mathbf
{I}_n)^{-1}\mathbf{x}_n\leq\frac{C}{\sqrt{v}}.
\]
Thus,
\begin{eqnarray*}
&&\mathrm{E}\bigl\vert m_{H^{\widehat{\mathbf{S}}_n}}(z)-m_{H^{\widetilde
{\mathbf{S}}_n}}(z)\bigr\vert
\\
&&\qquad \leq \mathrm{E} \llVert\widehat{\mathbf{S}}_n-\widetilde{
\mathbf{S}}_n\rrVert\bigl\llVert\mathbf{x}_n^* (
\widehat{\mathbf{S}}_n-z\mathbf{I}_n )^{-1}\bigr
\rrVert\bigl\llVert(\widetilde{\mathbf{S}}_n-z\mathbf{I}_n
)^{-1}\mathbf{x}_n\bigr\rrVert
\\
&&\qquad \leq C\sqrt{N}\eta_N^{-7} N^{-7/4} \mathrm{E}
\bigl( |X_{11}|^8 \mathrm{I} \bigl(|X_{11}|\geq
\eta_NN^{1/4} \bigr) \bigr)
\\
&&\quad\qquad{}\times \bigl(\mathrm{E} \bigl\|\mathbf{x}_n^*(\widehat{\mathbf{S}}_n-z
\mathbf{I})^{-1} \bigr\|^2 \bigr)^{1/2} \bigl( \mathrm{E}
\bigl\|(\widetilde{\mathbf{S}}_n-z\mathbf{I})^{-1}
\mathbf{x}_n \bigr\|^2 \bigr)^{1/2}
\\
&&\qquad \leq C\sqrt{N}v^{-3/2}\eta_N^{-7}
N^{-7/4} \mathrm{E} \bigl( |X_{11}|^8 \mathrm{I}
\bigl(|X_{11}|\geq\eta_NN^{1/4} \bigr) \bigr)
\\
&&\qquad \leq o\bigl(N^{-1/2}\bigr).
\end{eqnarray*}

\subsection{Rescaling}\label{sec9.4}
The rescaling procedures for the three theorems are exactly the same,
and only 8th moment is required. Thus,
we treat them uniformly. Write $\underline{\mathbf{Y}}_n=\widetilde
{\underline{\mathbf{X}}}_n/\sigma_1$,
where
\[
\sigma_1^2=\mathrm{E}\bigl\vert X_{11}
\mathrm{I} \bigl( |X_{11} |\leq\eta_NN^{1/4}
\bigr)-\mathrm{E} \bigl( X_{11}\mathrm{I} \bigl(|X_{11}|\leq
\eta_NN^{1/4} \bigr) \bigr)\bigr\vert^2.
\]
Notice that $\sigma_1$ tends to $1$ as $N$ goes to $\infty$. Define
$\mathbf{G}_n=\underline{\mathbf{Y}}_n\underline{\mathbf{Y}}_n^*/N$,
which is the sample covariance matrix of $\underline{\mathbf{Y}}_n$. We
shall show that $\mathbf{G}_n$ and $\mathbf{S}_n$ are asymptotically
equivalent, that is, the VESD of $\mathbf{G}_n$ and $\mathbf{S}_n$ have
the same limit if either one limit exists. For $v\geq C_0N^{-1/2}$,
\begin{eqnarray*}
\bigl\vert m_{H^{\mathbf{G}_n}}(z)-m_{H^{\widetilde{\mathbf
{S}}_n}}(z)\bigr\vert&=& \bigl
\vert\mathbf{x}_n^* (\widetilde{\mathbf{S}}_n-z\mathbf
{I}_n )^{-1} (\widetilde{\mathbf{S}}_n-
\mathbf{G}_n ) (\mathbf{G}_n-z\mathbf{I}_n
)^{-1}\mathbf{x}_n\bigr\vert
\\
&\leq& \frac{1}{v^2} \bigl\llVert\bigl(1-\sigma_1^{-1}
\bigr) \widetilde{\mathbf{S}}_n \bigr\rrVert\qquad(\mbox{see Lemma \ref{913}})
\\
&\leq& \frac{C}{v^2}\bigl(1-\sigma_1^2\bigr)\qquad\mbox{a.s.}
\\
&\leq& Cv^{-2} \eta_N^{-6} N^{-3/2}
\mathrm{E} \bigl(|X_{11}|^8 \mathrm{I}\bigl(|X_{11}|>
\eta_NN^{1/4}\bigr) \bigr)
\\
&\leq& o\bigl(N^{-1/2}\bigr)\qquad\mbox{a.s.}
\end{eqnarray*}

Hence, we shall without loss of generality assume that every $|X_{ij}|$
is bounded by $\eta_NN^{1/4}$, and every $X_{ij}$
has mean 0 and variance 1.
\end{appendix}



%

\printaddresses

\end{document}